# TESTING THE ORDER OF A MODEL

By Antoine Chambaz

*Université René Descartes*

This paper deals with order identification for nested models in the i.i.d. framework. We study the asymptotic efficiency of two generalized likelihood ratio tests of the order. They are based on two estimators which are proved to be strongly consistent. A version of Stein's lemma yields an optimal underestimation error exponent. The lemma also implies that the overestimation error exponent is necessarily trivial. Our tests admit nontrivial underestimation error exponents. The optimal underestimation error exponent is achieved in some situations. The overestimation error can decay exponentially with respect to a positive power of the number of observations.

These results are proved under mild assumptions by relating the underestimation (resp. overestimation) error to large (resp. moderate) deviations of the log-likelihood process. In particular, it is not necessary that the classical Cramér condition be satisfied; namely, the log-densities are not required to admit every exponential moment. Three benchmark examples with specific difficulties (location mixture of normal distributions, abrupt changes and various regressions) are detailed so as to illustrate the generality of our results.

**1. Introduction.** This paper is devoted to order identification problems in the independent and identically distributed (i.i.d.) framework. It fits in the general setting of model selection initiated by the seminal papers of Mallows [35], Akaike [1], Rissanen [38] and Schwarz [41]. Order identification deals with the estimation and test of a structural parameter which indexes the complexity of the common distribution of the observations. The purpose is to derive some new consistency and efficiency results. Order identification applies, for instance, to mixture models [42], where the order is (loosely speaking) the number of populations. Another example of application is abrupt changes models, where the order is (roughly) the number of changes.

It will be argued below that this example conveniently models a medical problem in which the order is the number of distinct levels of expression of a disease.

1.1. *Description of the problem.* We observe $n$ i.i.d. random variables $Z_1, \ldots, Z_n$ with values in a measurable sample space $(\mathcal{Z}, \mathcal{F})$ ($\mathcal{Z}$ is Polish). These observations are defined on a common measurable space upon which all the random variables will be defined.

The distribution $P^\star$ of $Z_1$ may belong to one model in the increasing family $\{\Pi_K\}_{K \geq 1}$ of nested models. Here, each $\Pi_K$ is a parametric collection of probability distributions which are absolutely continuous with respect to the same measure $\mu$,

$$\Pi_K = \{P_\theta : \theta \in \Theta_K\} \subset \Pi_{K+1},$$

where $\{(\Theta_K, d_K)\}_{K \geq 1}$ is an increasing family of nested metric parameter sets. In this paper $d_K$ will be abbreviated to $d$.

The integer $K$ is called the order of the model $\Pi_K$. It is also the order of any $P_\theta \in \Pi_K \setminus \Pi_{K-1}$ (with the convention $\Pi_0 = \varnothing$). The order of $P^\star$ is denoted by $K^\star$. It is infinite whenever $P^\star$ does not belong to $\Pi_\infty = \bigcup_{K \geq 1} \Pi_K$.

The central problem of this paper is an issue of composite hypotheses testing: we want to decide between the null hypothesis "$K^\star \leq K_0$" and its alternative "$K^\star > K_0$" (for some integer $K_0$), that is, to test

$$\text{"}P^\star \in \Pi_{K_0}\text{"} \quad \text{against} \quad \text{"}P^\star \notin \Pi_{K_0}\text{."}$$

This question is obviously crucial when the order is the quantity of interest. Furthermore, order identification may also be a prerequisite to consistent parameter estimation, when overestimation of the order causes loss of identifiability.

1.2. *Consistency and efficiency issues.* Let $\alpha_n$ and $\beta_n$ denote the type I and type II errors of a procedure that tests the hypotheses above. This procedure is consistent if $\alpha_n$ and $\beta_n$ converge to zero as $n$ tends to infinity. Its efficiency is measured in terms of rates of convergence of $\alpha_n$ and $\beta_n$ to zero.

In the classical statistical theory, a standard Neyman–Pearson procedure tests two simple hypotheses by comparing the log-likelihoods at each of them to a constant threshold. Now, it is known [10] that this procedure satisfies

$$\limsup_{n \to \infty} n^{-1} \log \alpha_n < 0 \quad \text{and} \quad \limsup_{n \to \infty} n^{-1} \log \beta_n < 0.$$

It is consequently natural, when investigating the efficiency of an order testing procedure, to study whether the rates of convergence are exponential with respect to $n$ or not.

Two generalized likelihood ratio test procedures based on two different estimators of $K^\star$ will be studied here. Obviously, if $\widetilde{K}_n$ estimates $K^\star$, then the natural rule is to reject the null hypothesis if $\widetilde{K}_n > K_0$. Then

$$\alpha_n \leq P^\star\{\widetilde{K}_n > K^\star\} \quad \text{and} \quad \beta_n \leq P^\star\{\widetilde{K}_n < K^\star\}$$

(these upper bounds do not depend on $K_0$). According to the discussion above, we shall thus focus on the following issues:

1. Are our order estimators strongly consistent?
2. Can $e_u > 0$ or $e_o > 0$ (the underestimation and overestimation error exponents, resp.) be found such that

$$\limsup_{n \to \infty} n^{-1} \log P^\star\{\widetilde{K}_n < K^\star\} \leq -e_u$$

or

$$\limsup_{n \to \infty} n^{-1} \log P^\star\{\widetilde{K}_n > K^\star\} \leq -e_o?$$

If so, can the error exponents $e_u$ or $e_o$ be arbitrarily large? If not, what happens at a subexponential rate, that is, when replacing the factor $n^{-1}$ by a factor $v_n^{-1} = o(1)$, with $v_n = o(n)$?

The consistency issue 1 has been studied for two decades. The interest in the efficiency issue 2 is more recent. By formulating the efficiency issue this way, we adopt the error exponent perspective of the information theory literature [13]. This notion of efficiency is asymptotic, as are all our results. It is connected to other notions of asymptotic efficiency, among which is Bahadur efficiency [3]. The latter is usually derived from large deviations results. In the following, the underestimation (resp. overestimation) error will similarly be related to large (resp. moderate) deviations of the log-likelihood process.

1.3. *Results in perspective.* Pioneering results about order identification of time series can be found in [2]. Strong consistency of the same order estimator in autoregressive models is shown in [24] and [26]. The test of the order of an ARMA process is addressed in [16]. Error exponents for autoregressive order testing are investigated in [8].

Consistent estimation of the order of a mixture model is at stake in [15, 21, 27, 29, 30, 34]. Efficiency issues are addressed in [15]. Also, [16] is concerned with the test of the order of a mixture.

Order estimation in exponential models is studied in [25]. The rates of underestimation and overestimation of two estimators of the order are investigated in [25] (for exponential models), [31] (for regular models) and in [23] (for models characterized by the existence of an exhaustive finite-dimensional statistic).

The problem of order identification in Markov models on a finite alphabet must be mentioned too. Some important papers are [12, 14] (they give insight into the consistency issue for some classical order estimators) and also [20, 22] (where optimal underestimation error exponents are obtained for the same classical order estimators). A more comprehensive presentation of order identification in Markov models can be found in [7].

*A new method for new results.* In most previous work the choice of the framework is contingent on the need for tractable explicit calculus. In this paper we shall resort to general properties of empirical processes. Our approach yields several new results that hold under mild assumptions.

In particular, our test procedures admit nontrivial underestimation error exponents. Besides, one of them has an optimal underestimation error exponent in some situations. Any test procedure based on a consistent estimator is proved to admit a necessarily trivial overestimation error exponent. The overestimation probabilities of our procedures can decay exponentially fast with respect to a positive power of $n$.

More details follow.

*Benchmark examples.* Let us introduce very briefly our three benchmark examples. Their presentation is merely sketched here, including the results obtained by applying our main general results. A whole section will be devoted to the detailed study of the examples.

Let $\sigma$ denote a known positive number.

- Location mixture example (LM): this is a notoriously difficult problem in the order identification literature (see the references cited above). In this model, one observes

$$Z_i = X_i + \sigma e_i \qquad (i = 1, \ldots, n),$$

where $X_1, \ldots, X_n$ are i.i.d. hidden (i.e., not observed) random variables with a common distribution of finite support $\{m_1, \ldots, m_{K^\star}\}$, and $e_1, \ldots, e_n$ are i.i.d. and independent from $X_1, \ldots, X_n$, with centered Gaussian distribution of variance 1. The goal is to estimate $K^\star$.

Applying the main general results of this paper will imply the following:

1. Our two estimators of $K^\star$ are consistent.
2. Their underestimation error exponents are nontrivial and bounded by a number which depends on squared distances between $P^\star$ and $\Pi_K$, $K = 1, \ldots, K^\star - 1$. Their overestimation error exponents are trivial but their overestimation probabilities decay exponentially fast with respect to a positive power of $n$.

   These results are new for maximum likelihood procedures.

- Abrupt changes example (AC): this example is original in the order identification literature. In this model one observes

$$Y_i = f^\star(X_i) + \sigma e_i \qquad (i = 1, \ldots, n),$$

where $X_1, \ldots, X_n$ are i.i.d. on a subset of $\mathbb{R}^q$ ($q \geq 2$); $e_1, \ldots, e_n$ are i.i.d. and independent of $X_1, \ldots, X_n$, with centered Gaussian distribution of variance 1, and the function $f^\star$ is piecewise constant. Loosely speaking, the goal is to estimate a minimal number of domains on which $f^\star$ is constant.

In virtue of the general results of this paper, the following new results hold ("almost surely" abbreviates to "a.s."):

1. $P^\star$-a.s., our estimators are greater than or equal to $K^\star$ eventually.
2. Our tests admit nontrivial underestimation error exponents. Their overestimation error exponents are necessarily trivial.

- Various regression examples (VR): let $\{t_k\}_{k \geq 1}$ be an orthonormal system in $L^2([0,1])$. In this model one observes

$$Y_i = f^\star(X_i) + \sigma e_i \qquad (i = 1, \ldots, n),$$

where $X_1, \ldots, X_n$ are i.i.d., uniformly distributed on $[0,1]$, $e_1, \ldots, e_n$ are i.i.d. and independent of $X_1, \ldots, X_n$, with centered Gaussian distribution of variance 1, and $f^\star = \sum_{k=1}^{K^\star} \theta_k t_k$ with $\theta_{K^\star} \neq 0$. The goal is to estimate $K^\star$.

As a consequence of the main general results of this paper, the following results are obtained:

1. Our two estimators of $K^\star$ are consistent.
2. Their underestimation error exponents are nontrivial, and one of them achieves optimality. Their overestimation error exponents are necessarily trivial, but their overestimation probabilities decay exponentially fast with respect to a positive power of $n$.

   In particular, the optimality of one of the underestimation error exponents is a new result.

1.4. *Organization of the paper.* In Section 2 some notation precedes the definition of the order estimators studied here. The basic assumptions are stated. Moreover, two limit theorems for the log-likelihood process which will play a central role are recalled. The consistency results are stated and commented on in Section 3. The most conclusive part is Section 4. It is devoted to the statement of the efficiency results and comments. The application of our general results to the benchmark examples is addressed in detail in Section 5. The proofs are postponed to the Appendix.

**2. Notation and preliminaries.** The integral $\int f \, d\lambda$ of a function $f$ with respect to a measure $\lambda$ will be written as $\lambda f$. Besides, all the expressions involving extrema and empirical processes will be assumed measurable.

2.1. *Two maximum penalized likelihood estimators.* Let $p_\theta$ denote the density of $P_\theta$ with respect to $\mu$ and $\ell_\theta = \log p_\theta$ (for all $\theta \in \Theta_\infty = \bigcup_{K \geq 1} \Theta_K$). $P^\star$ is supposed to be absolutely continuous with respect to $\mu$ without loss of generality. Its density is denoted by $p^\star$ and we set $\ell^\star = \log p^\star$. If $P^\star \in \Pi_{K^\star} \setminus \Pi_{K^\star - 1}$, then $P^\star = P_{\theta^\star}$ for $\theta^\star \in \Theta_{K^\star} \setminus \Theta_{K^\star - 1}$.

The log-likelihood $\ell_n$ of the observations is

$$\ell_n(\theta) = \sum_{i=1}^n \ell_\theta(Z_i) \qquad (\text{every } \theta \in \Theta_\infty).$$

The penalized maximum likelihood criterion for the model $\Pi_K$ is written as

$$\operatorname{crit}(n, K) = \sup_{\theta \in \Theta_K} \ell_n(\theta) - \operatorname{pen}(n, K),$$

where pen is a positive penalty function. It yields the two estimators of the order studied in this paper,

$$\widehat{K}_n^{\mathrm{L}} = \inf\{K \geq 1 : \operatorname{crit}(n, K) \geq \operatorname{crit}(n, K+1)\},$$
$$\widehat{K}_n^{\mathrm{G}} = \inf \arg\sup_{K \geq 1} \{\operatorname{crit}(n, K)\} \geq \widehat{K}_n^{\mathrm{L}}.$$

$\widehat{K}_n^{\mathrm{G}}$ is a global (hence, the G in its name) maximizer of the criterion. $\widehat{K}_n^{\mathrm{G}}$ always bounds from above $\widehat{K}_n^{\mathrm{L}}$, the first local (hence, the L) maximizer of the same criterion. Note that the computation of these estimators is a less demanding algorithmic task for $\widehat{K}_n^{\mathrm{L}}$ than for $\widehat{K}_n^{\mathrm{G}}$.

*Comment.* A prior bound $K_{\max}$ for $K^\star$ will be assumed known when studying the overestimation properties of $\widehat{K}_n^{\mathrm{G}}$. Indeed, we cannot control its overestimation probability when infinitely many models are involved. This assumption is common in the order identification literature [2, 7, 20, 21, 23, 24, 25, 30, 31].

On the one hand, there are situations where assuming the existence of $K_{\max}$ is mandatory. It is, for instance, proven [14] that some classical (minimum description length) order estimators are not consistent when no upper bound to the true order is known a priori: they fail to recover the true order 0 of a uniformly distributed i.i.d. sequence on a finite alphabet $A$, when $\Pi_K$ is the set of all Markov chains of order at most $K$. On the other hand, it is also shown in the same paper that the so-called Bayesian information criterion (BIC) order estimator is consistent when no upper bound is known a priori. It is thus particularly interesting that the study of the properties of $\widehat{K}_n^{\mathrm{L}}$ does not require a prior bound for $K^\star$.

Now, it must be emphasized that our asymptotic study of the problem does not allow us to obtain conditions on the dependence of $\operatorname{pen}(n, K)$ on $K$. In contrast, the former BIC order estimator studied by Csiszár and Shields [14] corresponds to $\operatorname{pen}(n, K) = \frac{1}{2}|A|^K(|A| - 1)\log n$. It is believed that this is a minimal penalty. In [22] the dependence on $K$ of the penalty function is also made precise (but the penalty is certainly not minimal, according to the authors).

The dependence of $\operatorname{pen}(n, K)$ on $K$ could be investigated through risk bounds for maximum log-likelihood [6, 36] (in the testing framework of this paper, the chosen loss function is $K \mapsto \mathbb{1}\{K \neq K^\star\}$). However, this would require at present time some restrictive assumptions. For instance, exact asymptotic risk bounds are yet out of reach for a mixture of Gaussian distributions. Furthermore, exact asymptotic bounds are not enough in overestimation, when we have to deal with infinitely many models [12].

2.2. *Basic assumptions.* Let us denote by $H(P|Q) = P \log dP/dQ$ if $P \ll Q$, $H(P|Q) = \infty$ otherwise, the relative entropy of $P$ with respect to $Q$. A survey of the relative entropy properties can be found, for instance, in [19]. If $\Pi$ is a subset of $M_1(\mathcal{Z})$ [the set of all the probability measures on $(\mathcal{Z}, \mathcal{F})$], the infimum of $H(P|Q)$ for $P$ (resp. $Q$) ranging through $\Pi$ will be denoted by $H(\Pi|Q)$ [resp. $H(P|\Pi)$].

The following assumptions will be needed throughout this paper:

A1. Compactness assumption. For all $K \geq 1$, the parameter sets $(\Theta_K, d)$ are compact metric sets and the models $\Pi_K$ are compact for the weak topology on the space $M_1(\mathcal{Z})$.
A2. Parameterization assumption. The parameterization $\theta \mapsto \ell_\theta(z)$ from $\Theta_K$ to $\mathbb{R}$ is continuous for all $z \in \mathcal{Z}$ and $K \geq 1$.
A3. Bracket assumption. There exist $l, u \in \mathbb{R}^{\mathcal{Z}}$ such that $(u - l) \in L^1(P^\star)$ and
$$l \leq \ell^\star \leq u \quad \text{and} \quad l \leq \ell_\theta \leq u \qquad (\text{all } \theta \in \Theta_\infty).$$
A4. Penalty assumption.
   $\operatorname{pen}(n, \cdot)$ is an increasing function for all $n \geq 1$.
   $\operatorname{pen}(n, K) \to \infty$ as $n \to \infty$ and $\operatorname{pen}(n, K) = o(n)$ for all $K \geq 1$.

The continuous parameterization assumption A2 is standard in statistics (see, e.g., [43]). Assumption A3 is called "bracket assumption" after the definition of the bracket $[l, u]$ (which is the set of all functions $f$ with $l \leq f \leq u$). It is also standard in the literature to invoke A3 when empirical processes are involved [43]. Another standard assumption in this setting is the boundedness of the parameter set. Assumption A1 is slightly stronger (at least when the parameter set is finite-dimensional, by virtue of the Heine–Borel theorem, A2 and Lévy's continuity theorem). Assumption A4 is the minimum requirement for a penalty function. Finally, it is worth noting that A3 implies that $H(P^\star|P_\theta)$ is finite for all $\theta \in \Theta_\infty$.

2.3. *Large and moderate deviation of the log-likelihood process.* It is shown in Section 4, which is devoted to efficiency issues, that underestimation can be related to large deviations of the log-likelihood process, while overestimation can be related to moderate deviations of the latter. Large and moderate deviations of the log-likelihood process both describe the limiting behavior of the empirical measure $\mathbb{P}_n = n^{-1} \sum_{i=1}^n \delta_{Z_i}$ ($\delta_z$ denotes the Dirac measure at $z$) on rare events as $n$ goes to infinity. Let us state the principles we shall need (their lower bounds are omitted).

*Extended Sanov theorem* [32]. Let $\tau$ be given by $\tau(s) = \exp(|s|) - |s| - 1$ (all $s \in \mathbb{R}$). The classes

(1) $\qquad \mathcal{L}_\tau(P^\star) = \{f \in \mathbb{R}^\mathcal{Z} : \exists a > 0, P^\star \tau(f/a) < \infty\},$

(2) $\qquad \mathcal{M}_\tau(P^\star) = \{f \in \mathbb{R}^\mathcal{Z} : \forall a > 0, P^\star \tau(f/a) < \infty\} \subset \mathcal{L}_\tau(P^\star)$

will play a central role in our study. $\mathcal{L}_\tau(P^\star)$ [resp. $\mathcal{M}_\tau(P^\star)$] is the set of all functions on $\mathcal{Z}$ that admit *some* (resp. *any*) exponential moment with respect to $P^\star$. In the LM example, for instance, if a continuous function $f$ upon $\mathbb{R}$ satisfies $f = O(x^2)$ at infinity, then $f \in \mathcal{L}_\tau(P^\star)$. For such a function, $f \in \mathcal{M}_\tau(P^\star)$ if and only if $f = o(x^2)$ at infinity. This simple example will be particularly interesting when $f$ is a log-density $\ell_\theta$ [which is an $O(x^2)$ but not an $o(x^2)$] or a difference $(\ell_\theta - \ell^\star)$ [which is an $o(x^2)$].

When equipped with the norm

(3) $\qquad \|f\|_\tau = \inf\{a > 0 : P^\star \tau(f/a) \leq 1\} \qquad (\text{all } f \in \mathcal{L}_\tau),$

$\mathcal{L}_\tau(P^\star)$ is a Banach space. Its topological dual is denoted by $\mathcal{L}'_\tau(P^\star)$. In this paper we shall be particularly interested in the set

$$\mathcal{Q} = \{Q \in \mathcal{L}'_\tau(P^\star) : Q \geq 0, Q1 = 1\} \cup \mathcal{P},$$

where $\mathcal{P} = \{p^{-1} \sum_{i=1}^p \delta_{z_i} : p \geq 1, z_1, \ldots, z_p \in \mathcal{Z}\}$. It is equipped with the coarsest topology that makes the linear forms $Q \mapsto Qf$ continuous for every $f \in \mathcal{L}_\tau(P^\star)$ and with the coarsest $\sigma$-field that makes them measurable. It is worth noting that $\mathbb{P}_n \in \mathcal{Q} \cap \mathcal{P} = \mathcal{P}$, hence, the need for $\mathcal{P}$.

By definition, $Q \in \mathcal{Q}$ is $P^\star$-*singular* if there exists a sequence $\{A_p\}$ of measurable sets such that $Q\mathbb{1}\{A_p^c\} = 0$ for all $p \geq 1$, while $\lim_{p \to \infty} P^\star(A_p) = 0$. It is known (Theorem 2.3 and Proposition 2.4 in [32]) that:

LEMMA 1. *Any $Q \in \mathcal{Q} \cap \mathcal{L}'_\tau(P^\star)$ is uniquely decomposed into the sum $Q = Q^a + Q^s$, where $Q^a \in \mathcal{L}'_\tau(P^\star)$ is a probability measure, $Q^a \ll P^\star$, while $Q^s \in \mathcal{L}'_\tau(P^\star)$ is $P^\star$-singular and $Q^s \geq 0$. Besides, for every $f \in \mathcal{M}_\tau(P^\star)$, $Qf = Q^a f$.*

REMARK 1. $\mathcal{Q} \cap \mathcal{L}'_\tau(P^\star)$ *is not* a subset of $M_1(\mathcal{Z})$. If $Q \in \mathcal{Q} \cap \mathcal{L}'_\tau(P^\star)$, then $P(A) = Q\mathbb{1}\{A\}$ (for any measurable set $A$) does define a probability measure $P$, which is in fact $Q^a$. Besides, $P$ and $Q$ coincide on $\mathcal{M}_\tau(P^\star)$, but may differ on $\mathcal{L}_\tau(P^\star) \setminus \mathcal{M}_\tau(P^\star)$ ($Q = P = Q^a$ if and only if $Q^s = 0$).

Let us finally introduce the nonnegative function $I$ (the *extended* relative entropy) defined for any $Q = Q^a + Q^s \in \mathcal{Q} \cap \mathcal{L}'_\tau(P^\star)$ by

$$I(Q) = H(Q^a|P^\star) + \sup\{Q^s f : f \in \mathcal{L}_\tau(P^\star), P^\star \exp(f) < \infty\}$$

and $I(Q) = \infty$ if $Q \in \mathcal{Q} \cap \mathcal{P} = \mathcal{P}$. It particularly satisfies the following:

LEMMA 2. *For every $Q \in \mathcal{Q}$, $I(Q) \geq 0$, with equality if and only if $Q = P^\star$.*

Theorem 3.2 in [32] encompasses the following result.

THEOREM 1 [32]. *The function $I$ is a convex, lower semicontinuous mapping from $\mathcal{Q}$ to $[0,\infty]$. Its level sets $\{Q \in \mathcal{Q} : I(Q) \leq \alpha\}$ are compact for all $\alpha > 0$. Moreover, for any measurable $S \subset \mathcal{Q}$ [with closure $\mathrm{cl}(S)$],*

$$\limsup_{n \to \infty} n^{-1} \log P^\star\{\mathbb{P}_n \in S\} \leq - \inf_{Q \in \mathrm{cl}(S)} I(Q).$$

REMARK 2. Theorem 1 requires an involved setting. Three reasons motivate its use, though:

- A classical Sanov theorem on $M_1(\mathcal{Z})$ would be insufficient here. Indeed, when dealing with the underestimation rate, our proofs require that the linear forms $Q \mapsto Q\ell_\theta$ be continuous on $\mathcal{Q}$ (any $\theta \in \Theta_\infty$), while possibly $\ell_\theta \in \mathcal{L}_\tau(P^\star) \setminus \mathcal{M}_\tau(P^\star)$. Now, Schied [40] has shown that the extension of a Sanov theorem on $M_1(\mathcal{Z})$ to a topology on $M_1(\mathcal{Z})$ that makes the linear form $Q \mapsto Qf$ continuous on $\mathcal{Q}$ for some $f \in \mathcal{L}_\tau(P^\star)$ is possible if and only if $f \in \mathcal{M}_\tau(P^\star)$ (this is the classical Cramér condition).
- Provided the need that $Q \mapsto Qf$ be continuous on $\mathcal{Q}$ for various $f \in \mathcal{L}_\tau(P^\star) \setminus \mathcal{M}_\tau(P^\star)$, the topology on $\mathcal{Q}$ introduced above is the natural one.
- The simpler relative entropy rate function $I'(Q) = H(Q^a|P^\star)$ for $Q = Q^a + Q^s \in \mathcal{Q} \cap \mathcal{L}'_\tau(P^\star)$, $I'(Q) = \infty$ otherwise, does not have compact level sets (this is also a consequence of [40]). This would be a major drawback in our scheme of proof.

*Moderate deviations of* $\mathbb{P}_n$ [44]. Let $\mathcal{G}$ denote a subclass of $L^2(P^\star)$ with envelope $G \in \mathbb{R}^{\mathcal{Z}}$ [i.e., $|g(z)| \leq G(z)$ for all $g \in \mathcal{G}$ and $z \in \mathcal{Z}$].

Let $\ell^\infty(\mathcal{G})$ be the collection of all bounded functions $b \in \mathbb{R}^{\mathcal{G}}$. The uniform norm $\|\cdot\|_{\mathcal{G}}$ defined by $\|b\|_{\mathcal{G}} = \sup_{g \in \mathcal{G}} |b(g)|$ induces a topology and a $\sigma$-field on $\ell^\infty(\mathcal{G})$.

Let us denote by $M_0(\mathcal{Z})$ the space of all signed measures $Q$ on $(\mathcal{Z}, \mathcal{F})$ that satisfy $Q1 = 0$, $\sup_{g \in \mathcal{G}} |Qg| < \infty$ and $Q \ll P^\star$ (the derivative $dQ/dP^\star$ is denoted by $q$). One observes that, for any $Q \in M_0(\mathcal{Z})$, $Q^\infty g = Qg$ (all $g \in \mathcal{G}$) defines an element of $\ell^\infty(\mathcal{G})$. Particularly, $(\mathbb{P}_n - P^\star)^\infty$ is a random variable on $\ell^\infty(\mathcal{G})$ under $P^\star$.

Let us finally introduce the nonnegative function $J$ defined for any $b \in \ell^\infty(\mathcal{G})$ by

$$J(b) = \inf\left\{P^\star \frac{q^2}{2} : Q \in M_0(\mathcal{Z}), Q^\infty = b\right\}$$

(with the convention $\inf \varnothing = +\infty$).

THEOREM 2 [44]. *Let $\{v_n\}$ be an increasing sequence of positive numbers such that $v_n = o(n)$, $n \log n = o(v_n^2)$. Let us assume that there exist $A \geq 1$, $\delta \in (0,1)$ such that, for every $k, n \geq 1$,*

$$v_{nk} \leq A k^{1-\delta} v_n.$$

*If $\mathcal{G}$ is $P^\star$-Donsker and $G \in \mathcal{L}_\tau(P^\star)$, then for any $S \subset (\ell^\infty(\mathcal{G}), \|\cdot\|_{\mathcal{G}})$,*

$$\limsup_{n \to \infty} (v_n^2/n)^{-1} \log P^\star \{nv_n^{-1}(\mathbb{P}_n - P^\star)^\infty \in S\} \leq -\inf_{b \in \mathrm{cl}(S)} J(b).$$

This theorem is a straightforward corollary of Theorem 5 in [44] (for a recent account of the $P^\star$-Donsker property, see [43]).

**3. Consistency issue.** The statements of our three results of consistency are gathered here. These results are rather routine. However, the resort to empirical process arguments allows us to achieve great generality. We refer to Section 5 for examples of application and comparison with previous consistency results in each benchmark framework.

From now on, Log denotes the truncated log, that is, $\mathrm{Log}(x) = \log(x \vee e)$ (all $x \in \mathbb{R}$). The function $\varphi$ is defined by $\varphi(x) = x^2/\mathrm{Log}\,\mathrm{Log}(x)$ (all $x \in \mathbb{R}$). Besides, let us introduce the classes of functions

(4) $\qquad \mathcal{G}_K^{\mathrm{a}} = \{g_\theta = (\ell_\theta - \ell^\star) : \theta \in \Theta_K\} \qquad$ (every $K \geq 1$).

THEOREM 3. *Let $P^\star$ belong to $\Pi_{K^\star} \setminus \Pi_{K^\star-1}$. Suppose that $\varphi(u-l) \in L^1(P^\star)$ and that the penalty function satisfies*

$$\liminf_{n\to\infty} \frac{\mathrm{pen}(n, K+1)}{\mathrm{pen}(n, K)} > 1 \quad \text{and} \quad \limsup_{n\to\infty} \frac{(n \log \log n)^{1/2}}{\mathrm{pen}(n, K)} = 0 \qquad (\text{any } K \geq 1).$$

- *If $P^\star \notin \Pi_K$ implies $H(P^\star | \Pi_{K+1}) < H(P^\star | \Pi_K)$, and if, moreover, $\mathcal{G}^{\mathrm{a}}_{K^\star+1}$ is $P^\star$-Donsker, then $P^\star$-a.s., $\widehat{K}^{\mathrm{L}}_n = K^\star$ eventually.*
- *If $K^\star \leq K_{\max}$, and if, moreover, $\mathcal{G}^{\mathrm{a}}_{K_{\max}}$ is $P^\star$-Donsker, then $P^\star$-a.s., $\widehat{K}^{\mathrm{G}}_n = K^\star$ eventually.*

It is proved in Section 5 that the theorem applies to the LM and VR examples. In the AC example, it is obtained that $P^\star$-a.s., $\widehat{K}^{\mathrm{L}}_n \geq K^\star$ and $\widehat{K}^{\mathrm{G}}_n \geq K^\star$ eventually (see Proposition B.1).

The scheme of proof of the latter theorem is rather standard. The proof of "no underestimation eventually" relies on the strong law of large numbers. It essentially requires the continuous parameterization assumption A2 and finally boils down to a comparison of the following:

- $H(P^\star | \Pi_K)$ with $H(P^\star | \Pi_{K+1})$ for all $K < K^\star - 1$ when dealing with $\widehat{K}^{\mathrm{L}}_n$ (hence, the assumption that strict inequality holds);
- $H(P^\star | \Pi_{K^\star-1}) > 0$ with $H(P^\star | \Pi_{K^\star}) = 0$ when dealing with $\widehat{K}^{\mathrm{G}}_n$ (this comparison is obvious).

The proof of "no overestimation eventually" relies on a law of the iterated logarithm. It essentially requires the $P^\star$-Donsker assumptions.

We emphasize that the condition on the penalty function in Theorem 3 excludes BIC-like expressions $\mathrm{pen}(n, K) = \frac{1}{2} \dim(\Theta_K) \log n$. This can be overcome, as shown in Theorem 4, by resorting to an example of a "peeling device" (see Appendix A). To this end, substitutes for $\mathcal{G}^{\mathrm{a}}_K$ classes are introduced, namely,

$$(5) \qquad \mathcal{G}^{\mathrm{b}}_K = \left\{ g_\theta = \frac{\ell_\theta - \ell^\star}{H(\theta)^{1/2}} : \theta \in \Theta_K, H(\theta) > 0 \right\} \qquad (\text{every } K \geq 1),$$

where $H(\theta) = H(P^\star | P_\theta)$ for all $\theta \in \Theta_\infty$.

THEOREM 4. *Let $P^\star$ belong to $\Pi_{K^\star} \setminus \Pi_{K^\star-1}$. Suppose that $\varphi(u - l) \in L^1(P^\star)$ and that the penalty function satisfies*

$$\liminf_{n\to\infty} \frac{\mathrm{pen}(n, K+1)}{\mathrm{pen}(n, K)} > 1 \quad \text{and} \quad \limsup_{n\to\infty} \frac{\log \log n}{\mathrm{pen}(n, K)} = 0 \qquad (\text{any } K \geq 1).$$

- *If $P^\star \notin \Pi_K$ implies $H(P^\star | \Pi_{K+1}) < H(P^\star | \Pi_K)$, and if, moreover, $\mathcal{G}^{\mathrm{b}}_{K^\star+1}$ is $P^\star$-Donsker, then $P^\star$-a.s., $\widehat{K}^{\mathrm{L}}_n = K^\star$ eventually.*

- If $K^\star \leq K_{\max}$, and if, moreover, $\mathcal{G}^{\mathrm{b}}_{K_{\max}}$ is $P^\star$-Donsker, then $P^\star$-a.s., $\widehat{K}^{\mathrm{G}}_n = K^\star$ eventually.

This theorem applies to the LM example as proved in Section 5.
Finally, the last result of this section addresses a case of misspecification.

THEOREM 5. *Suppose that, for every $K \geq 1$, $0 < H(P^\star|\Pi_{K+1}) < H(P^\star|\Pi_K)$. Then $P^\star \notin \Pi_\infty$ and $P^\star$-a.s.,*

$$\liminf_{n \to \infty} \widehat{K}^{\mathrm{L}}_n = \liminf_{n \to \infty} \widehat{K}^{\mathrm{G}}_n = \infty.$$

The proofs are postponed to Appendix B.

**4. Efficiency issues.** It is argued in Section 1 that the efficiency issue in order identification problems is related to the decay to zero of $P^\star\{\widetilde{K}_n < K^\star\}$ and $P^\star\{\widetilde{K}_n > K^\star\}$ (for $\widetilde{K}_n = \widehat{K}^{\mathrm{L}}_n$ or $\widehat{K}^{\mathrm{G}}_n$) as $n$ tends to infinity. A comparison with standard Neyman–Pearson tests suggested investigating whether they can vanish exponentially fast with respect to $n$ or not. This is the question at stake in the next section.

4.1. *Best error exponents.* We shall resort hereafter to a concise version of Stein's lemma (our Lemma 3) due to Bahadur, Zabell and Gupta [5] (see their Theorem 2.1, specialized here to the case of i.i.d. processes for sake of simplicity). An early version is mentioned in [10] in a framework of hypotheses testing, and stated, for instance, in [4]. Lemma 3 relies on the core of Stein's original proof (which is a change of probability argument). It is, in most cases, the key of its various versions.

LEMMA 3 [5]. *Let $\mathbb{P}, \mathbb{Q}$ be two probability measures on the same measured space and $\{X_n\}$ be a sequence of random variables on it. Let $\{A_n\}$ be a sequence of measurable sets such that $A_n$ is $\sigma(X_1, \ldots, X_n)$-measurable.*

*Assume that $\mathbb{P} = P^{\otimes \infty}$ and $\mathbb{Q} = Q^{\otimes \infty}$, so that $\{X_n\}$ is an i.i.d. process under $\mathbb{P}$ and $\mathbb{Q}$.*

(6) $\quad$ *If* $\liminf_{n \to \infty} \mathbb{Q}(A_n) > 0, \quad$ *then* $\liminf_{n \to \infty} n^{-1} \log \mathbb{P}(A_n) \geq -H(Q|P).$

Now, by virtue of Lemma 3:

THEOREM 6. *Let $\widetilde{K}_n$ be any estimator of the order of the common distribution of $Z_1, \ldots, Z_n$.*

*Underestimation. If for all $K_0 \geq 1$ and $P_0 \in \Pi_{K_0} \setminus \Pi_{K_0 - 1}$,*

(7) $$\limsup_{n \to \infty} P_0\{\widetilde{K}_n > K_0\} < 1,$$

*then*

$$\liminf_{n\to\infty} n^{-1}\log P^\star\{\widetilde{K}_n < K^\star\} \geq -\inf_{K<K^\star} H(\Pi_K|P^\star) = -H(\Pi_{K^\star-1}|P^\star).$$

Overestimation. *If for all $K_0 \geq 1$ and $P_0 \in \Pi_{K_0} \setminus \Pi_{K_0-1}$,*

$$\limsup_{n\to\infty} P_0\{\widetilde{K}_n < K_0\} < 1,$$

*then*

$$\liminf_{n\to\infty} n^{-1}\log P^\star\{\widetilde{K}_n > K^\star\} = \limsup_{n\to\infty} n^{-1}\log P^\star\{\widetilde{K}_n > K^\star\} = 0.$$

PROOF. Set $K_0 < K^\star$ and $\theta_0 \in \Theta_{K_0}$. Choose the probabilities $P = P^\star$, $Q = P_{\theta_0}$ and define $A_n = (\widetilde{K}_n \leq K_0)$.

The left-hand side condition of (6) is satisfied by virtue of (7), hence the right-hand side property of (6) holds. Now, $P^\star\{\widetilde{K}_n < K^\star\} \geq P^\star\{A_n\}$ and $K_0, \theta_0$ are arbitrary, so the proof in the underestimation case is complete.

The proof in the overestimation case parallels the lines above. □

Analogous versions of this theorem have been proved in [20] and [22] in settings of Markov chains and hidden Markov models order identification, respectively. It is, however, and surprisingly a new result (to the best of our knowledge) in our framework of order identification from i.i.d. observations. In summary, the underestimation (resp. overestimation) result holds for estimators $\widetilde{K}_n$ that ultimately overestimate (resp. underestimate) the order with a probability bounded away from one. Thus, the theorem applies to any consistent estimator. Besides, the conclusion of Theorem 6 for such estimators is twofold:

- The underestimation probability can decay exponentially fast with respect to $n$, and a best possible underestimation error exponent, namely, $H(\Pi_{K^\star-1}|P^\star)$, is exhibited.
- The overestimation probability cannot decay exponentially fast with respect to $n$: the overestimation error exponent is necessarily trivial.

Consequently, the main issue is now to prove that $\widehat{K}_n^{\mathrm{L}}$ and $\widehat{K}_n^{\mathrm{G}}$ admit nontrivial underestimation error exponents and to compare those exponents to $H(\Pi_{K^\star-1}|P^\star)$. This will involve large deviations of the log-likelihood process; see Section 4.2. The second issue is to investigate the behavior of the overestimation probabilities. These probabilities will be related to moderate (instead of large) deviations of the latter process; see Section 4.3.

4.2. *Underestimation error exponent.* Let us introduce for any $\alpha \geq 0$ and $K \geq 1$ the following subsets of $\mathcal{Q}$ ($\Lambda$ stands for Local and $\Gamma$ for Global):

$$\Lambda_{\alpha,K} = \left\{ Q \in \mathcal{Q} : \sup_{\theta \in \Theta_K} Q\ell_\theta - \sup_{\theta \in \Theta_{K+1}} Q\ell_\theta \geq -\alpha \right\}, \tag{8}$$

$$\Gamma_{\alpha,K} = \left\{ Q \in \mathcal{Q} : \sup_{\theta \in \Theta_K} Q\ell_\theta - \sup_{\theta \in \Theta_{K^\star}} Q\ell_\theta \geq -\alpha \right\}. \tag{9}$$

$\{\Gamma_{\alpha,K}\}$ is nondecreasing in $\alpha$ and $K$, and $\{\Lambda_{\alpha,K}\}$ is nondecreasing in $\alpha$. Besides, for every $\alpha \geq 0$ and $K < K^\star$, $\Gamma_{\alpha,K} \subset \Lambda_{\alpha,K}$. Finally, $\Gamma_{\alpha,K^\star} = \mathcal{Q}$.

From now on, let us suppose that $P^\star \in \Pi_{K^\star} \setminus \Pi_{K^\star - 1}$.

THEOREM 7. *Let us assume that, for every $\theta \in \Theta_\infty$, $H(P_\theta | P^\star)$ is finite and $p_\theta \ell_\theta \in L^1(\mu)$. Let us also suppose that:*

(i) $\{\ell_\theta : \theta \in \Theta_\infty\} \subset \mathcal{L}_\tau(P^\star)$.

(ii) *For all $K \geq 1$, for every $Q \in \mathcal{Q}$ and $\varepsilon > 0$ small enough, there exists a finite subset $T \subset \Theta_K$ such that*

$$\forall \theta \in \Theta_K, \exists t \in T : |Q\ell_\theta - Q\ell_t| \leq \varepsilon.$$

- *If $P^\star \notin \Pi_K$ implies $H(P^\star | \Pi_{K+1}) < H(P^\star | \Pi_K)$, then*

$$\limsup_{n \to \infty} n^{-1} \log P^\star \{ \widehat{K}_n^{\mathrm{L}} < K^\star \} \leq - \inf_{K < K^\star} I(\Lambda_{0,K}) < 0. \tag{10}$$

- *Moreover,*

$$\limsup_{n \to \infty} n^{-1} \log P^\star \{ \widehat{K}_n^{\mathrm{G}} < K^\star \} \leq -I(\Gamma_{0,K^\star - 1}) < 0. \tag{11}$$

*Furthermore, if $P^\star$-a.s. for any $n \geq 1$, $Z_1, \ldots, Z_n$ are mutually distinct, then "for every $Q \in \mathcal{Q}$" may be replaced by "for every $Q \in \mathcal{Q} \cap \mathcal{L}'_\tau(P^\star)$" in (ii) [ yielding (ii)$'$].*

This theorem fully applies to the LM, AC and VR examples, as proved in Section 5.

REMARK 3. The alternative assumption (ii)$'$ is needed for the AC example. The proof of the theorem is slightly more involved with the relaxed condition. It particularly requires a more precise framework for the large deviations principle of Theorem 1 (refer to the proof in Section C.1 for further details).

*Comment on Theorem* 7. Theorem 7 is the most conclusive result of this paper. It notably relates the phenomenon of underestimation to the large deviations of the log-likelihood process. The assumptions of the theorem are mild and give, we think, insight into the phenomenon of underestimation. This assertion is justified by the fact that the assumptions are satisfied in the three benchmark examples, despite their differences and specific difficulties. This is due to the resort to empirical processes arguments (and recent advances in large deviations theory) in place of tractable explicit calculus. Let us emphasize that Theorem 7 applies even when the log-densities $\ell_\theta$ admit *some* exponential moment rather than *any* (the classical Cramér condition).

Comparison with previous results on the rate of underestimation in each benchmark framework can be found in Section 5.

Besides, a comparison with [23] is relevant. In the latter, the authors consider an order estimator based on the minimization over a *finite-dimensional* parameter set of an empirical criterion $U_n(\theta)$. The basic assumption requires the existence of a *finite-dimensional* statistic $T_n$ which satisfies an *exponential maximal inequality* and the existence of a continuous function $U$ such that $U_n(\theta) = U(\theta, T_n)$. In this framework, tractable calculus in finite dimensions yields some nonasymptotic evaluation of the underestimation probability (and the overestimation probability too). The scope of the paper is large, although its basic assumption excludes mixture models (and particularly LM) because there is no finite-dimensional statistic $T_n$ for the log-likelihood; it also excludes models with infinite-dimensional parameter sets (and particularly AC).

*Optimal underestimation error exponent.* We aim at showing that, under appropriate further assumptions, $\widehat{K}_n^{\mathrm{G}}$ achieves the optimal underestimation error exponent. Theorem 8 is an intermediate result, in which *possibly* tighter upper bounds for the probabilities of underestimation are stated.

Let us reinforce the structure of the spaces $\Theta_K$: now, the distance $d$ on $\Theta_K$ derives from a norm $\|\cdot\|$ on the vector space $\Theta_K$, so that the notion of differentiability with respect to $\theta \in \Theta_K$ is available. This particularly excludes the AC example.

THEOREM 8. *Suppose that the assumptions of Theorem* 7 *are valid. In addition, assume that:*

(iii) $(u - l) \in \mathcal{M}_\tau(P^\star)$.

(iv) *For every $K \leq K^\star$ and $z \in \mathcal{Z}$, the functions $\theta \mapsto \ell_\theta(z)$ are differentiable on the interior of $\Theta_K$, with derivative $\dot{\ell}_\theta(z)$. Moreover, the coordinates of $\dot{\ell}_\theta$ are elements of $\mathcal{M}_\tau(P^\star)$ and there exists $F \in \mathcal{L}_\tau(P^\star)$ such that*

$$(12) \qquad |\ell_{\theta+h} - \ell_\theta - \dot{\ell}_\theta^T h| \leq F \cdot o(h).$$

- If $P^\star \notin \Pi_K$ implies $H(P^\star|\Pi_{K+1}) < H(P^\star|\Pi_K)$, then

$$\limsup_{n\to\infty} n^{-1} \log P^\star\{\widehat{K}_n^{\mathrm{L}} < K^\star\} \leq - \inf_{K<K^\star} H(\Lambda_{0,K} \cap M_1(\mathcal{Z})|P^\star) < 0.$$

- Moreover,

$$\limsup_{n\to\infty} n^{-1} \log P^\star\{\widehat{K}_n^{\mathrm{G}} < K^\star\} \leq -H(\Gamma_{0,K^\star-1} \cap M_1(\mathcal{Z})|P^\star) < 0.$$

This theorem fully applies in the LM and VR examples, as proved in Section 5.

REMARK 4. In assumption (iv), inequality (12) may hold only for $h = \|h\|e_k$, where $e_k$ is the $k$th canonical basis vector of $\Theta_K$.

In conclusion, the underestimation error exponent turns out to be optimal (regarding Theorem 6) for $\widehat{K}_n^{\mathrm{G}}$ in exponential models. This is another new result. It applies particularly to our sole exponential model, that is in the VR example, as shown in Section 5.

THEOREM 9. *Under the assumptions of Theorem* 8 *and* for exponential models, *the best underestimation error exponent (regarding Theorem* 6*) is achieved by* $\widehat{K}_n^{\mathrm{G}}$.

*Comment on Theorems* 8 *and* 9. It is easily seen that the upper bounds of Theorem 8 are indeed lower than the ones in Theorem 7, but possibly *not strictly*. Are there situations where these inequalities are known to be strict or not strict? What is the nature of the discrepancy between the optimal exponent and the one obtained in Theorem 7? These are very difficult questions, to which we do not have any answer. Boucheron and Gassiat [8] faced the same impediment when they studied the underestimation error exponent of a procedure which tests the order of an autoregressive process. They first show that their order estimator has nontrivial underestimation error exponent. A version of Stein's lemma yields an optimal error exponent. They finally check that, in some situations, the optimal error exponent is achieved. In both the present work and theirs, the main difficulty stems from the absence of a "full information-theoretical interpretation" (we quote their expression) of the large deviations rate function—that is, stems from the discrepancy between the rate function and the relative entropy.

It is also worth emphasizing that, although the optimal underestimation efficiency is proved for $\widehat{K}_n^{\mathrm{G}}$ in exponential models (see Theorem 9), we cannot conclude that $\widehat{K}_n^{\mathrm{L}}$ and $\widehat{K}_n^{\mathrm{G}}$ are not optimal on the basis of Theorems 7 and 8. In greater generality, we are not aware of any example of order estimator proven to be suboptimal regarding the underestimation error exponent in the statistical or information theoretical literatures.

The proof of Theorem 9 (postponed to Section C.3) involves *H-projections* as defined and studied by Csiszár [11]: $\overline{Q}$ is the $H$-projection of the probability measure $Q$ on a convex set of probability measures $\mathcal{C}$ [which must satisfy $H(P|Q) < \infty$ for some $P \in \mathcal{C}$] if $\overline{Q} \in \mathcal{C}$ and

(13) $$H(\overline{Q}|Q) = H(\mathcal{C}|Q).$$

$H$-projections satisfy a useful characterization (see Theorem 2.2 in [11]):

LEMMA 4 [11]. $Q' \in \mathcal{C}$ with $H(Q'|Q) < \infty$ *is the $H$-projection of $Q$ on $\mathcal{C}$ if and only if, for every $P \in \mathcal{C}$,*

$$H(P|Q) \geq H(P|Q') + H(Q'|Q).$$

Nonetheless, the proof also involves probability measures $P, \overline{P}$ and a set $\mathcal{C}$ which satisfy $\overline{P} \in \mathcal{C}$ and

(14) $$H(P|\overline{P}) = H(P|\mathcal{C}).$$

We shall say by analogy that $\overline{P}$ is the reversed-$H$-projection of $P$ on $\mathcal{C}$. Such reversed-$H$-projections are much less tractable than $H$-projections. In general, notably, a reversed-$H$-projection cannot be characterized as in Lemma 4. However, it is remarkable that, *in exponential models*, reversed-$H$-projections do satisfy a similar characterization (the proof draws its inspiration from [9]):

LEMMA 5. *Set $Q \in \mathcal{L}'_\tau(P^\star) \cap M_1(\mathcal{Z})$ such that $H(Q|\Pi_{K^\star}) < \infty$. Then $Q \ll \mu$ (let $q$ denote its density $\mathrm{d}Q/\mathrm{d}\mu$). Now, let us assume that:*

(i) $\Pi_{K^\star}$ *is an exponential model.*
(ii) $\{\ell_\theta : \theta \in \Theta_{K^\star}\} \subset \mathcal{L}_\tau(P^\star)$.
(iii) $Q \log q < \infty$.
(iv) *The function $\theta \mapsto Q\ell_\theta$ is continuous from $\Theta_{K^\star}$ to $\mathbb{R}$.*

*Let $P$ belong to $\Pi_{K^\star}$. $P$ is the reversed-$H$-projection of $Q$ on $\Pi_{K^\star}$ if and only if*

(15) $$H(Q|P_\theta) \geq H(Q|P) + H(P|P_\theta) \qquad (any\ P_\theta \in \Pi_{K^\star}).$$

These characterizations will play a central role in the proof of Theorem 9.

4.3. *Overestimation rate.* The following theorem provides a first link between the penalization function and the rate of overestimation (which is necessarily slower than exponential in $n$; see Theorem 6) that it yields for $\widehat{K}_n^{\mathrm{L}}$ and $\widehat{K}_n^{\mathrm{G}}$.

THEOREM 10. *Let the penalty function be of the form* $\mathrm{pen}(n,K) = v_n D(K)$, *where* $D \in \mathbb{R}^{\mathbb{N}}$ *and* $\{v_n\}$ *increase,* $v_n = o(n)$, *and for some* $A \geq 1$, $\delta \in (0,1)$, *for every* $k, n \geq 1$,

$$v_{nk} \leq A k^{1-\delta} v_n.$$

*Let us also suppose that:*

(i) $(u - l) \in \mathcal{L}_\tau(P^\star)$, *so that the classes* $\mathcal{G}_K^{\mathrm{a}}$ *[defined in* (4)*] admit an envelope function in* $\mathcal{L}_\tau(P^\star)$.
(ii) $n = o(v_n^2)$.

- *If* $\mathcal{G}_{K^\star+1}^{\mathrm{a}}$ *is* $P^\star$-*Donsker, then*

(16) $$\limsup_{n \to \infty} n v_n^{-2} \log P^\star \{\widehat{K}_n^{\mathrm{L}} > K^\star\} < 0.$$

- *If* $K^\star \leq K_{\max}$, *and if, moreover,* $\mathcal{G}_{K_{\max}}^{\mathrm{a}}$ *is* $P^\star$-*Donsker, then*

(17) $$\limsup_{n \to \infty} n v_n^{-2} \log P^\star \{\widehat{K}_n^{\mathrm{G}} > K^\star\} < 0.$$

For instance, $v_n = n^{1-\delta}$, $\delta \in (0, 1/2)$ is an admissible sequence and Theorem 10 applies to the LM and VR example.

The resort to the same "peeling device" that allowed the transition from Theorem 3 to Theorem 4 (both devoted to the consistency issue) in Section 3 yields again a relaxed condition on $\{v_n\}$.

THEOREM 11. *Let* pen *be of the form detailed in Theorem* 10. *Let us also suppose that:*

(i) *The classes* $\mathcal{G}_K^{\mathrm{b}}$ *[defined in* (5)*] admit an envelope function in* $\mathcal{L}_\tau(P^\star)$.
(ii) $\log n = o(v_n)$.

- *If* $\mathcal{G}_{K^\star+1}^{\mathrm{b}}$ *is* $P^\star$-*Donsker, then*

(18) $$\limsup_{n \to \infty} v_n^{-1} \log P^\star \{\widehat{K}_n^{\mathrm{L}} > K^\star\} < 0.$$

- *If* $K^\star \leq K_{\max}$, *and if, moreover,* $\mathcal{G}_{K_{\max}}^{\mathrm{b}}$ *is* $P^\star$-*Donsker, then*

(19) $$\limsup_{n \to \infty} v_n^{-1} \log P^\star \{\widehat{K}_n^{\mathrm{G}} > K^\star\} < 0.$$

For instance, $v_n = (\log n)^{1+\epsilon}$ ($\epsilon > 0$) is an admissible sequence.

*Comment on Theorems* 10 *and* 11. Theorem 10 is the main result on the efficiency issue of overestimation in this paper. It notably relates the phenomenon of overestimation to the moderate deviations of the log-likelihood process. The assumptions of the theorem are rather mild. This opinion is justified by the fact that the theorem applies to the LM and VR examples. It is worth pointing out that the conditions related to the log-densities $\ell_\theta$ are expressed in terms of the envelope function and the $P^\star$-Donsker property (and not in terms of exponential moments for $\ell_\theta$). As explained in Section 5, the AC example is excluded because we do not verify the $P^\star$-Donsker property.

On the contrary, Theorem 11 relies on strong assumptions, particularly assumption (i), which exclude the LM and VR examples. Although the condition on $\{v_n\}$ is relaxed, Theorem 11 does not apply to the BIC-like penalty function $\mathrm{pen}(n,K) = \frac{1}{2}\dim(\Theta_K)\log n$ ($v_n = \log n$). Besides, it is important to note that the choice of $v_n = (\log n)^{1+\varepsilon}$ yields control of the overestimation probability that decays like a negative power of $n$.

We refer to Section 5 for comparison with previous results on the rate of overestimation in the LM and VR benchmark examples (none exists for the AC example). The last paragraph of the comment of Theorem 7 is also relevant here, as a paradigm of the methods based on tractable calculus in finite dimensions.

**5. Benchmark examples.** This section is devoted to a detailed investigation of our benchmark examples in order to illustrate the collection of results that have been stated in the two previous sections.

5.1. *Location mixture example.* Let $\sigma$ be a priori known and $\gamma(\cdot;m)$ denote the density of the Gaussian distribution with mean $m$ and variance $\sigma^2$ with respect to the Lebesgue measure $\mu$ on $\mathbb{R}$. Let $\mathcal{M}$ be a compact subset of $\mathbb{R}$. Here, $\Pi_1$ is the set of all Gaussian probability measures with mean $m \in \mathcal{M}$ and variance $\sigma^2$ and $\Theta_1 = \mathcal{M}$. For every $\theta \in \Theta_1$, let us define $p_\theta = \gamma(\cdot;\theta)$. Now, for any $K \geq 2$, let us introduce the compact sets

$$\Theta_K = \left\{\theta = (\boldsymbol{\pi}, \mathbf{m}) : \boldsymbol{\pi} = (\pi_1, \ldots, \pi_{K-1}) \in \mathbb{R}_+^{K-1}, \sum_{k=1}^{K-1} \pi_k \leq 1, \mathbf{m} \in \mathcal{M}^K \right\}.$$

Every $\theta \in \Theta_K$ ($K \geq 2$) is associated with a mixing distribution $F_\theta = \sum_{k=1}^{K-1} \pi_k \times \delta_{m_k} + (1 - \sum_{k=1}^{K-1})\delta_{m_K}$ on $\mathcal{M}$ and a probability measure $P_\theta$ with density $p_\theta = \int_{\mathcal{M}} \gamma(\cdot;m)\,dF_\theta(m)$ with respect to $\mu$. For $K \geq 2$, $\Pi_K = \{P_\theta : \theta \in \Theta_K\}$.

In this setting, one observes

$$Z_i = X_i + \sigma e_i \qquad (i=1,\ldots,n),$$

where $X_1,\ldots,X_n$ are i.i.d. hidden random variables, $e_1,\ldots,e_n$ are i.i.d. and independent of $X_1,\ldots,X_n$, with centered Gaussian distribution of variance 1, and there exists $\theta^\star \in \Theta_{K^\star} \setminus \Theta_{K^\star-1}$ such that $X_1,\ldots,X_n$ have distribution $F_{\theta^\star}$. In this case, $Z_1,\ldots,Z_n$ are i.i.d. and $P^\star$-distributed.

*Exploring the assumptions.* The compactness assumption A1 is easily verified (by virtue of Lévy's continuity theorem). The continuous parameterization assumption A2 is satisfied. Defining $l = \inf \ell_\theta$ and $u = \sup \ell_\theta$ (the suprema range over $\theta \in \Theta_\infty$) ensures $l \leq \ell_\theta \leq u$ (all $\theta \in \Theta_\infty$) and $(u-l)^{1+c} \in \mathcal{L}_\tau(P^\star)$ for some $c > 0$. Hence, the bracket assumption A3 holds. Now, a slight adaptation of the proof of Lemma 3 in [34] yields the following:

PROPOSITION 1. *Let $F$ be a mixing distribution on $\mathcal{M}$ (possibly with infinite support) and $P^\star$ have density $p^\star = \int_\mathcal{M} \gamma(\cdot; m) \, dF(m)$. In the LM example, if $P^\star \notin \Pi_K$, then $H(P^\star|\Pi_{K+1}) < H(P^\star|\Pi_K)$.*

The classes $\mathcal{G}_K^a$ are $P^\star$-Donsker (indeed, Example 19.7 in [43] guarantees that they have finite bracketing entropy integral). It can also be proven by hand that the classes $\mathcal{G}_K^b$ are $P^\star$-Donsker too (they have $\varepsilon$-bracketing numbers bounded by a polynomial in $\varepsilon^{-1}$, hence, finite bracketing entropy integral; see [43] for details). Consequently, the consistency conclusions of Theorems 3 and 4 are valid.

As for the efficiency issue of the underestimation rate, the assumptions of Theorems 7 and 8 are verified in this example. If $P^\star \in \Pi_{K^\star} \setminus \Pi_{K^\star-1}$, it is clear that, for every $\theta \in \Theta_\infty$, $H(P_\theta|P^\star)$ is finite, $p_\theta \ell_\theta \in L^1(\mu)$ and $\ell_\theta \in \mathcal{L}_\tau(P^\star)$ [this is assumption (i) of Theorem 7]. Moreover, as proved in Section E.1 (essentially by virtue of Ascoli's theorem applied to the restrictions of the $\ell_\theta$'s to a compact set) we have the following:

LEMMA 6. *In the LM example, the finite sieve assumption* (ii) *of Theorem 7 is satisfied.*

Assumption (iii) of Theorem 8 holds because $(u-l)^{1+c} \in \mathcal{L}_\tau(P^\star)$, hence, $(u-l) \in \mathcal{M}_\tau(P^\star)$. Furthermore, it can be shown (resorting to Taylor's integral remainder formula, e.g.) that assumption (iv) of Theorem 8 also holds, so that the latter applies in the LM example.

Finally, the assumptions of Theorem 10, which deal with the efficiency issue of the overestimation rate, have already been verified above.

In summary, Theorems 3, 4, 6, 7, 8 and 10 apply in the LM example.

*Comment.* Order identification in mixture models, even with known standard deviation, is a notoriously difficult problem.

Mixture models have been postulated in many applications; see Chapter 2 of [42] for a scope of these applications. Mixture models are notably characterized by their lack of identifiability when overestimating the order, and the subsequent singularity of the Fisher information matrix, which prevents one from using classical methods based on a Taylor expansion. They are also known for the tediousness of the related calculus. Besides, the log-densities $\ell_\theta$ *do not* belong to $\mathcal{M}_\tau(P^\star)$ (the strong Cramér condition is not satisfied), hence, the need for theorems that apply to the case of $\ell_\theta \in \mathcal{L}_\tau(P^\star) \setminus \mathcal{M}_\tau(P^\star)$.

Let us review here the previous results of order identification in mixture models that can be found in the literature. Regarding the consistency issue, Henna [27], Dacunha-Castelle and Gassiat [15] and James, Priebe and Marchette [29] proved the consistency (without any prior bound on the true order) of three different order estimators which do not rely on a maximum likelihood procedure. The consistency of our estimator $\widehat{K}_n^{\mathrm{G}}$ (with prior bound) has been already proven in this setting in [30] and [21]. The proof in [30] involves the locally conic parameterization of a mixture model introduced in [16] (this parameterization allows one to cope with Taylor expansions). The proof of [21] relies on a clever inequality for likelihood ratios which makes her proof very simple.

As for the efficiency issue, we have made clear in Sections 4.2 and 4.3 that large or moderate deviations of the log-likelihood process are a reasonable (and certainly minimal) requirement in order to yield asymptotic bounds on the probabilities of underestimation and overestimation. Hence, the locally conic parameterization does not appear adequate to yield such bounds. Dacunha-Castelle and Gassiat [15] proved that their estimator $\widetilde{K}_n = \arg\max_K \{U_n(K) + \mathrm{pen}(n, K)\}$ [where $U_n(K)$ depends on the data and $K$ and differs from the log-likelihood maximized on $\Theta_K$, using our notation of Section 4.3 for pen] satisfies, for some $c_1, c_2 > 0$ and $n$ large enough, $P^\star\{\widetilde{K}_n \neq K^\star\} \leq c_1 \exp(-c_2 n^{-1} v_n^2)$. The corresponding rate is the one of Theorem 10.

As far as we know, our results on efficiency stated in Theorems 6, 7, 8 and 10 are new for our maximum likelihood procedures.

5.2. *Abrupt changes example.* Let $(\mathcal{X}, \mathcal{B}, P)$ be an open subset of $\mathbb{R}^q$ ($q \geq 2$) equipped with the trace $\mathcal{B}$ of the Borel $\sigma$-field and a probability measure $P \ll \mu$, the Lebesgue measure on $\mathcal{X}$ (with density $\mathrm{d}P/\mathrm{d}\mu$ denoted by $p$).

Let CP be the set of all countable Caccioppoli partitions of $\mathcal{X}$. It is known that there exists a metric $d$ on CP such that the subset $\mathrm{CP}_b$ of all partitions whose "perimeters" are bounded by a fixed constant $b > 0$ is a compact metric space when equipped with $d$. (The definitions and main properties of Caccioppoli partitions can be found in [33].)

A partition is a family $\tau = \{\tau_j\}_{j \geq 1}$ of measurable subsets of $\mathcal{X}$ such that $P(\mathcal{X} \setminus \bigcup_j \tau_j) = 0$, $P(\tau_j \cap \tau_{j'}) = 0$ for every $j \neq j'$ and possibly $P(\tau_j) = 0$. The cardinality of $\tau$ is the number of $j \geq 1$ such that $P(\tau_j) > 0$. Given a compact set $\mathcal{M}$ of $\mathbb{R}$ and $\tau \in \mathrm{CP}_b$, it is easy to verify that one can associate $m_j \in \mathcal{M}$ with every $\tau_j$, yielding a marked partition $\{(\tau_j, m_j)\}_{j \geq 1}$, then modify the definition of $d$ so that the set of all marked partitions of $\mathrm{CP}_b$ is also a compact set when equipped with $d$. It is worth noting that, if $d[(\tau^0, m^0), (\tau^1, m^1)] \leq \delta$, then there exists a bijective map $\varphi$ from $I^0 = \{j : P(\tau_j^0) > 0\}$ to $\{j : P(\tau_j^1) > 0\}$ such that $P(\tau_j^0 \Delta \tau_{\varphi(j)}^1) \leq \delta$ and $|m_j^0 - m_{\varphi(j)}^1| \leq \delta$ for every $j \in I^0$ ($\Delta$ denotes the symmetrical difference between sets).

In this example, for every $K \geq 1$, $\Theta_K$ is the set of all marked partitions of $\mathrm{CP}_b$ with cardinality at most $K$. $(\Theta_K, d)$ is a compact metric space, hence, the first half of the compactness assumption A1. For $\sigma$ a priori known, let us denote by $\gamma(\cdot; m)$ the density of the Gaussian distribution with mean $m$ and variance $\sigma^2$, $f_\theta(x) = \sum_{k \geq 1} m_k \mathbb{1}\{x \in \tau_k\}$ and finally $p_\theta(z) = \gamma(y; f_\theta(x)) p(x)$ [for all $z = (x, y) \in \mathcal{Z} = \mathcal{X} \times \mathbb{R}$, $K \geq 1$ and $\theta \in \Theta_K$]. Let $P_\theta$ have $p_\theta$ for density with respect to $\mu$, then set $\Pi_K = \{P_\theta : \theta \in \Theta_K\}$ for every $K \geq 1$.

In this setting, one observes $Z_i = (X_i, Y_i)$ with

$$Y_i = f^\star(X_i) + \sigma e_i \qquad (i = 1, \ldots, n),$$

where $X_1, \ldots, X_n$ are i.i.d. and $P$-distributed, $e_1, \ldots, e_n$ are i.i.d. and independent of $X_1, \ldots, X_n$, with centered Gaussian distribution of variance 1, and there exists $\theta^\star \in \Theta_{K^\star} \setminus \Theta_{K^\star - 1}$ such that $f^\star = f_{\theta^\star}$. In this case, $Z_1, \ldots, Z_n$ are i.i.d. and $P_{\theta^\star}$-distributed.

*Exploring the assumptions.* Lévy's continuity theorem implies that the second half of A1 is satisfied. Besides, the continuous parameterization assumption A2 is obviously verified. It is easily seen that the bracket assumption A3 holds. Indeed, if one introduces $\underline{f} = \inf \ell_\theta$ and $\overline{f} = \sup \ell_\theta$ (the suprema range over $\Theta_\infty$), functions $l, u \in \mathbb{R}^{\mathcal{Z}}$ can be defined such that $(u - l)$ is continuous, $l \leq \ell_\theta \leq u$ (all $\theta \in \Theta_\infty$) and $2\sigma^2 (u - l)(z) = (\underline{f}^2 + \overline{f}^2)(x) + 2|y|(\overline{f} - \underline{f})(x)$, hence, $(u - l)^{1+c} \in \mathcal{L}_\tau(P^\star)$ for some $c > 0$.

Furthermore, if the $L^2(P)$-norm is denoted by $\|\cdot\|_2$, then it is worth stressing that, for every $\theta, t \in \Theta_\infty$,

$$(20) \qquad H(P_\theta | P_t) = \frac{\|f_\theta - f_t\|_2^2}{2\sigma^2}.$$

Using (20) yields (the proof is postponed to Section E.2) the following:

LEMMA 7. *In the AC example, if $P^\star \in \Pi_\infty \setminus \Pi_K$, then $H(P^\star | \Pi_{K+1}) < H(P^\star | \Pi_K)$.*

At this stage, Proposition B.1 of Section B applies. The proposition guarantees that underestimation eventually does not occur almost surely. On the contrary, Proposition B.2 does not apply because the required $P^\star$-Donsker properties are not verified. Thus, the overestimation probability cannot be controlled.

As for the efficiency issue of the underestimation rate, the assumptions of Theorem 7 are valid in this example. If $P^\star \in \Pi_{K^\star} \setminus \Pi_{K^\star-1}$, it is clear that, for every $\theta \in \Theta_\infty$, $H(P_\theta|P^\star)$ is finite, $p_\theta \ell_\theta \in L^1(\mu)$ and $\ell_\theta \in \mathcal{L}_\tau(P^\star)$ [this is assumption (i) of Theorem 7]. Finally, for all $n \geq 1$, $Z_1, \ldots, Z_n$ are mutually different $P^\star$-a.s. and, as shown in Section E.3,

LEMMA 8. *In the AC example, the finite sieve assumption* (ii)′ *of Theorem 7 is satisfied.*

In summary, $P^\star$-a.s., $\widehat{K}_n^{\mathrm{G}} \geq \widehat{K}_n^{\mathrm{L}} \geq K^\star$ eventually and the part of Theorem 6 which deals with overestimation and Theorem 7 apply in the AC example.

*Comment.* This example is original in the order identification literature. It is related to variational image segmentation theory, although it does not entirely fit in this general framework (because we observe the random responses $Y_i$ at random points $X_i$ rather than the random responses $Y_x$ at all $x \in \mathcal{X}$). This framework of order identification is a priori difficult, notably because the parameter sets $\Theta_K$ are not finite-dimensional.

The following medical problem is conveniently modeled by the AC example. Let us suppose that a disease is characterized by distinct levels of expression $k = 1, \ldots, K^\star$, whose number $K^\star$ is unknown. Let us also assume that:

- The mean of a clinical measure $Y$ (modeled by a Gaussian random variable of known variance $\sigma^2$) is uniquely characterized by the level $k$ of expression of the disease.
- Simultaneously, there exist $q \geq 2$ feature (demographic, diet, clinical) measurements $(x^1, \ldots, x^q) \in \mathcal{X}$ and a segmentation $\tau^\star = (\tau_k^\star)_{1 \leq k \leq K^\star}$ of the space $\mathcal{X}$ of their possible values, so that each $\tau_k^\star$ corresponds uniquely to the level $k$ of the disease.

Then, if one observes both $X_i = (X_i^1, \ldots, X_i^q)$ and $Y_i$ for $i = 1, \ldots, n$ patients, one may wish to estimate the number $K^\star$ of distinct levels of the disease.

5.3. *Various regressions example.* Let $\{t_K\}_{K \geq 1}$ be a uniformly bounded system of continuous functions on $[0, 1]$. Let us also assume that it is an orthonormal system in $L^2([0, 1])$ (equipped with Lebesgue measure). Let $\sigma$ be a priori known and $\gamma(\cdot; m)$ be the density of the Gaussian distribution with mean $m$ and variance $\sigma^2$. Let $\mathcal{M}$ be a compact subset of $\mathbb{R}$ that contains 0. Let us define $\Theta_K = \mathcal{M}^K$ (each $K \geq 1$). For every $\theta \in \Theta_K$, let us set $f_\theta = \sum_{k=1}^K \theta_k t_k$ and $p_\theta(z) = \gamma(y; f_\theta(x))$ (all $z = (x, y) \in [0, 1] \times \mathbb{R}$). Let $P_\theta$ have $p_\theta$ for density with respect to Lebesgue measure on $[0, 1] \times \mathbb{R}$, then set $\Pi_K = \{P_\theta : \theta \in \Theta_K\}$.

In this setting, one observes $Z_i = (X_i, Y_i)$ with

$$Y_i = f^\star(X_i) + \sigma e_i \qquad (i = 1, \ldots, n),$$

where $X_1, \ldots, X_n$ are i.i.d. and uniformly distributed on $[0,1]$, $e_1, \ldots, e_n$ are i.i.d. and independent of $X_1, \ldots, X_n$, with centered Gaussian distribution of variance 1, and there exists $\theta^\star \in \Theta_{K^\star} \setminus \Theta_{K^\star-1}$ such that $f^\star = f_{\theta^\star}$. In this case, $Z_1, \ldots, Z_n$ are i.i.d. and $P^\star$-distributed.

*Exploring the assumptions.* The compactness assumption A1 is clearly satisfied (by virtue of Lévy's continuity theorem for $\Pi_K$). Besides, the continuous parameterization assumption A2 is readily verified. The bracket assumption A3 holds: with $\underline{f} = \inf \ell_\theta$ and $\overline{f} = \sup \ell_\theta$ (the suprema range over $\theta \in \Theta_\infty$), $l, u \in \mathbb{R}^{\mathcal{Z}}$ can be defined such that $(u - l)$ is continuous, $l \leq \ell_\theta \leq u$ (any $\theta \in \Theta_\infty$) and $2\sigma^2(u-l)(z) = (\underline{f}^2 + \overline{f}^2)(x) + 2|y|(\overline{f} - \underline{f})(x)$, hence, $(u-l)^{1+c} \in \mathcal{L}_\tau(P^\star)$ for some $c > 0$. We emphasize that equality (20) also holds in this example when $\|\cdot\|_2$ denotes the $L^2([0,1])$ norm. A straightforward consequence follows:

LEMMA 9. *In the VR example, if $P^\star \in \Pi_\infty \setminus \Pi_K$, then $H(P^\star|\Pi_{K+1}) < H(P^\star|\Pi_K)$.*

Now, it can be proven that the classes $\mathcal{G}_K^a$ (all $K \geq 1$) defined in (4) are $P^\star$-Donsker (by mimicking the proof in the LM example), hence, the consistency conclusions of Theorem 3 are valid.

As for the efficiency issue of the underestimation rate, the assumptions of Theorems 7, 8 and 9 are satisfied in this example. If $P^\star \in \Pi_{K^\star} \setminus \Pi_{K^\star-1}$, it is clear that, for every $\theta \in \Theta_\infty$, $H(P_\theta|P^\star)$ is finite, $p_\theta \ell_\theta \in L^1(\mu)$ and $\ell_\theta \in \mathcal{L}_\tau(P^\star)$ [this is assumption (i) of Theorem 7]. Moreover, following the proof of Lemma 6 yields the following:

LEMMA 10. *In the VR example, the finite sieve assumption* (ii) *of Theorem 7 is satisfied.*

Now, it has been already argued that $(u-l)^{1+c} \in \mathcal{L}_\tau(P^\star)$, hence, $(u-l) \in \mathcal{M}_\tau(P^\star)$ and assumption (iii) of Theorem 8 is valid. Furthermore, a crude yet careful application of Taylor's integral remainder theorem yields assumption (iv) of Theorem 8. In conclusion, the models are exponential in the VR example, so Theorem 9 applies.

Concerning the efficiency issue of the overestimation rate, the assumptions of Theorem 10 have been verified in the lines above.

In summary, Theorems 3, 6, 7, 8, 9 and 10 apply in the VR example.

*Comment.* We present this example because it fits in the general framework of order identification in nested exponential models. This important framework has been investigated in [25] and [31] (who actually address the more general case of regular models). In the latter, the authors study the properties of $\widehat{K}_n^{\mathrm{G}}$ (with a prior bound on the true order). They prove its weak consistency. Rates of underestimation and overestimation similar to the ones of Theorems 7 and 10 are obtained. However, the underestimation error exponent is not shown to be at most $H(\Pi_{K^\star-1}|P^\star)$ and, of course, is not compared to it.

Thus, to the best of our knowledge, the results of Theorem 3, 6, 7, 8 and 10 are new in this exponential model framework for $\widehat{K}_n^{\mathrm{L}}$ (which does not require any prior bound on the true order), while the results of Theorems 3, 6 and 9 are new for $\widehat{K}_n^{\mathrm{G}}$.

## APPENDIX A: AN EXAMPLE OF THE PEELING DEVICE

The so-called "peeling device" classically allows one to analyze the rate of convergence of $M$-estimators in nonclassical frameworks. The original idea is due to Huber [28]. Examples may be found, for instance, in [37] for simple proofs of uniform central limit theorems or in [6] (see Proposition 7 therein and the attached remark) in a framework of risk bounds model selection. Another form of this device is the core of [21], where it applies to an order estimation problem for a mixture with Markov regime.

PROPOSITION A.1. *Set $K_2 > K_1 \geq K^\star$, the order of $P^\star$. Then, both inequalities below hold, the second one providing an example of the peeling technique:*

$$(\text{A.1}) \qquad \sup_{\theta \in \Theta_{K_2}} |(\mathbb{P}_n - P^\star)(\ell_\theta - \ell^\star)| \geq \sup_{\theta \in \Theta_{K_2}} \mathbb{P}_n \ell_\theta - \sup_{\theta \in \Theta_{K_1}} \mathbb{P}_n \ell_\theta$$

*and*

$$(\text{A.2}) \qquad \left( \sup_{\theta \in \Theta_{K_2}} \left| (\mathbb{P}_n - P^\star) \frac{\ell_\theta - \ell^\star}{H(\theta)^{1/2}} \right| \right)^2 \geq \sup_{\theta \in \Theta_{K_2}} \mathbb{P}_n \ell_\theta - \sup_{\theta \in \Theta_{K_1}} \mathbb{P}_n \ell_\theta.$$

PROOF. Inequality (A.1) is readily proved, since

$$\sup_{\theta \in \Theta_{K_2}} \mathbb{P}_n \ell_\theta - \sup_{\theta \in \Theta_{K_1}} \mathbb{P}_n \ell_\theta \leq \sup_{\theta \in \Theta_{K_2}} \mathbb{P}_n (\ell_\theta - \ell^\star)$$
$$= \sup_{\theta \in \Theta_{K_2}} \{(\mathbb{P}_n - P^\star)(\ell_\theta - \ell^\star) + P^\star(\ell_\theta - \ell^\star)\}$$
$$\leq \sup_{\theta \in \Theta_{K_2}} (\mathbb{P}_n - P^\star)(\ell_\theta - \ell^\star).$$

For (A.2), let us define for all $\theta \in \Theta_{K_2}$ such that $H(\theta) > 0$ (i.e., $P^\star \neq P_\theta$) the scaled log-densities ratio

$$g_\theta = \frac{\ell_\theta - \ell^\star}{H(\theta)^{1/2}}$$

and $g_\theta = 0$ otherwise. Now, for any $\theta \in \Theta_{K_2}$, $H(\theta)$ nonnegative yields

(A.3)
$$\mathbb{P}_n(\ell_\theta - \ell^\star) + H(\theta) = (\mathbb{P}_n - P^\star)(\ell_\theta - \ell^\star)$$
$$\leq H(\theta)^{1/2} \sup_{\theta \in \Theta_{K_2}} (\mathbb{P}_n - P^\star)g_\theta.$$

Let us set some $\theta_0 \in \Theta_{K_2}$ such that both $\sup_{\theta \in \Theta_{K_2}} \mathbb{P}_n(\ell_\theta - \ell^\star) \leq \mathbb{P}_n(\ell_{\theta_0} - \ell^\star) + \varepsilon$ and $\mathbb{P}_n(\ell_{\theta_0} - \ell^\star) \geq 0$. Then, (A.3) implies, for $\theta = \theta_0$,

$$\sup_{\theta \in \Theta_{K_2}} \mathbb{P}_n(\ell_\theta - \ell^\star) \leq H(\theta_0)^{1/2} \sup_{\theta \in \Theta_{K_2}} (\mathbb{P}_n - P^\star)g_\theta + \varepsilon.$$

Furthermore, $\mathbb{P}_n(\ell_{\theta_0} - \ell^\star) \geq 0$ combined with (A.3) imply in turn

$$H(\theta_0) \leq H(\theta_0)^{1/2} \sup_{\theta \in \Theta_{K_2}} (\mathbb{P}_n - P^\star)g_\theta,$$

hence,

$$\sup_{\theta \in \Theta_{K_2}} \mathbb{P}_n \ell_\theta - \sup_{\theta \in \Theta_{K_1}} \mathbb{P}_n \ell_\theta \leq \sup_{\theta \in \Theta_{K_2}} \mathbb{P}_n(\ell_\theta - \ell^\star) \leq \left( \sup_{\theta \in \Theta_{K_2}} (\mathbb{P}_n - P^\star)g_\theta \right)^2 + \varepsilon,$$

which completes the proof, since $\varepsilon > 0$ is arbitrary. $\square$

## APPENDIX B: PROOFS OF CONSISTENCY

**B.1. No underestimation eventually.** A strong law of large numbers for the supremum of the likelihood ratios is stated. Its routine proof relies on the achievement of $H(P^\star|\Pi_K)$, the standard strong law of large numbers and the Borel–Lebesgue property.

LEMMA B.1. $P^\star$-a.s., for any $K \geq 1$,
$$\sup_{\theta \in \Theta_K} n^{-1}(\ell_n(\theta) - \ell_n(\theta^\star)) \xrightarrow[n \to \infty]{} -H(P^\star|\Pi_K).$$

Now the result of no underestimation can be stated and proved. It is seen in Section 5 that Proposition B.1 fully applies to the LM and AC examples. It is also shown that the VR example satisfies the assumption in the case of $\widehat{K}_n^{\mathrm{G}}$.

PROPOSITION B.1. *Let us assume that $P^\star \in \Pi_{K^\star} \setminus \Pi_{K^\star - 1}$.*

- If $P^\star \notin \Pi_K$ implies $H(P^\star|\Pi_{K+1}) < H(P^\star|\Pi_K)$, then $P^\star$-a.s., $\widehat{K}_n^{\mathrm{L}} \geq K^\star$ eventually.
- $P^\star$-a.s., $\widehat{K}_n^{\mathrm{G}} \geq K^\star$ eventually.

PROOF. Let us abbreviate "infinitely often" to i.o. and prove that $P^\star\{\widehat{K}_n^{\mathrm{L}} < K^\star \text{ i.o.}\} = 0$ (minor changes allow us to cope with $\widehat{K}_n^{\mathrm{G}}$). By the union bound, it suffices to show that $P^\star\{\widehat{K}_n^{\mathrm{L}} = K \text{ i.o.}\} = 0$ for $K = 1, \ldots, K^\star - 1$. Now, if we denote by $\delta = H(P^\star|\Pi_{K+1}) - H(P^\star|\Pi_K) < 0$,

$$P^\star\{\widehat{K}_n^{\mathrm{L}} = K \text{ i.o.}\} \leq P^\star\left\{\sup_{\theta \in \Theta_K} \mathbb{P}_n \ell_\theta - \sup_{\theta \in \Theta_{K+1}} \mathbb{P}_n \ell_\theta \geq -\delta/2 \text{ i.o.}\right\}$$

$$\leq P^\star\left\{\liminf_{n\to\infty}\left\{\sup_{\theta \in \Theta_K} \mathbb{P}_n \ell_\theta - \sup_{\theta \in \Theta_{K+1}} \mathbb{P}_n \ell_\theta\right\} \geq -\delta/2\right\},$$

where the first inequality stems from the definition of the penalty function A4 and is satisfied for $n$ large enough. Finally, Lemma B.1 ensures that the right-hand side probability is zero, which concludes the proof. $\square$

The proof of Theorem 5 also fits in this "no underestimation" section.

PROOF OF THEOREM 5. $P^\star \notin \Pi_\infty$ because otherwise there would exist a $K \geq 1$ such that $H(P^\star|\Pi_K) = 0$. Lemma B.1 implies that, $P^\star$-a.s. and for all $K \geq 1$,

$$\sup_{\theta \in \Theta_{K+1}} \mathbb{P}_n \ell_\theta - \sup_{\theta \in \Theta_K} \mathbb{P}_n \ell_\theta \xrightarrow[n\to\infty]{} H(P^\star|\Pi_K) - H(P^\star|\Pi_{K+1}) > 0.$$

Therefore, by virtue of the definition of the penalty function A4, $P^\star$-a.s.,

$$\mathrm{crit}(n, K+1) - \mathrm{crit}(n, K) \xrightarrow[n\to\infty]{} \infty,$$

hence, $\widehat{K}_n^{\mathrm{G}} \geq \widehat{K}_n^{\mathrm{L}} > K$ for $n$ large enough. This is true for any $K \geq 1$, so the proof is complete. $\square$

**B.2. No overestimation eventually.**

PROPOSITION B.2. *Let us assume that $P^\star \in \Pi_{K^\star} \setminus \Pi_{K^\star - 1}$.*

- *If $\varphi(u-l) \in L^1(P^\star)$, if $\mathcal{G}^{\mathrm{a}}_{K^\star+1}$ (resp. $\mathcal{G}^{\mathrm{b}}_{K^\star+1}$) is $P^\star$-Donsker, then whenever* pen *satisfies the condition of Theorem 3 (resp. Theorem 4), $P^\star$-a.s., $\widehat{K}_n^{\mathrm{L}} \leq K^\star$ eventually.*
- *Let $K^\star \leq K_{\max}$. If $\varphi(u-l) \in L^1(P^\star)$, if $\mathcal{G}^{\mathrm{a}}_{K_{\max}}$ (resp. $\mathcal{G}^{\mathrm{b}}_{K_{\max}}$) is $P^\star$-Donsker, then whenever* pen *satisfies the condition of Theorem 3 (resp. Theorem 4), $P^\star$-a.s., $\widehat{K}_n^{\mathrm{G}} \leq K^\star$ eventually.*

It is proven in Section 5 that the assumptions of the proposition are satisfied in both cases a (i.e., under the assumptions of Theorem 3) and b (i.e., under the assumptions of Theorem 4) in the LM example. In the VR example, they are satisfied in case a.

The following lemma is a bounded law of the iterated logarithm stated in convenient terms for our purpose. It is a simple consequence of Theorem 4.1. in [18]. It is involved in the proof of Proposition B.2.

LEMMA B.2 [18]. *Let us assume that $\varphi(u - l) \in L^1(P^\star)$ and that, for some $K > K^\star$, $\mathcal{G} = \mathcal{G}_K^{\mathrm{a}}$ (resp. $\mathcal{G} = \mathcal{G}_K^{\mathrm{b}}$) is $P^\star$-Donsker. Then there exists a positive constant $C_K$ such that, $P^\star$-a.s.,*

$$\limsup_{n \to \infty} \frac{n^{1/2} \sup_{g \in \mathcal{G}} |(\mathbb{P}_n - P^\star)g|}{(\log \log n)^{1/2}} \leq C_K.$$

PROOF OF PROPOSITION B.2. Set $K = K^\star + 1$.

$P^\star\{\widehat{K}_n^{\mathrm{L}} > K^\star \text{ i.o.}\}$

$$\leq P^\star\bigg\{\sup_{\theta \in \Theta_K} \mathbb{P}_n \ell_\theta - \sup_{\theta \in \Theta_{K^\star}} \mathbb{P}_n \ell_\theta \geq n^{-1}\{\mathrm{pen}(n,K) - \mathrm{pen}(n,K^\star)\} \text{ i.o.}\bigg\}$$

$$\leq P^\star\bigg\{\bigg[\frac{(n \log \log n)^{1/2}}{\mathrm{pen}(n,K^\star)}\bigg]\bigg[\frac{n^{1/2} \sup_{g \in \mathcal{G}_K^{\mathrm{a}}} |(\mathbb{P}_n - P^\star)g|}{(\log \log n)^{1/2}}\bigg]$$

$$\geq \frac{\mathrm{pen}(n,K)}{\mathrm{pen}(n,K^\star)} - 1 \text{ i.o.}\bigg\},$$

where the last inequality is straightforward [it is (A.1)]. Consequently, whenever $\varphi(u - l) \in L^1(P^\star)$ and $\mathcal{G}_K^{\mathrm{a}}$ is $P^\star$-Donsker, Lemma B.2 applies and implies that, if pen satisfies the condition of Theorem 3, then $P^\star\{\widehat{K}_n^{\mathrm{L}} > K^\star \text{ i.o.}\} = 0$.

Now, renormalization yields an alternative bound for the second probability in the display above [by using (A.2) of the peeling technique Proposition A.1], namely

$P^\star\{\widehat{K}_n^{\mathrm{L}} > K^\star \text{ i.o.}\}$

$$\leq P^\star\bigg\{\bigg[\frac{\log \log n}{\mathrm{pen}(n,K^\star)}\bigg]\bigg[\frac{n^{1/2} \sup_{g \in \mathcal{G}_K^{\mathrm{b}}} |(\mathbb{P}_n - P^\star)g|}{(\log \log n)^{1/2}}\bigg]^2 \geq \frac{\mathrm{pen}(n,K)}{\mathrm{pen}(n,K^\star)} - 1 \text{ i.o.}\bigg\}.$$

Therefore, if $\varphi(u - l) \in L^1(P^\star)$ and $\mathcal{G}_K^{\mathrm{b}}$ is $P^\star$-Donsker, Lemma B.2 applies and implies that, as soon as pen satisfies the condition of Theorem 3, $P^\star\{\widehat{K}_n^{\mathrm{L}} > K^\star \text{ i.o.}\} = 0$. This concludes the study of $\widehat{K}_n^{\mathrm{L}}$.

Furthermore, if $K^\star \leq K_{\max}$, then the union bound guarantees that it suffices to prove that $P^\star\{\widehat{K}_n^{\mathrm{G}} = K \text{ i.o.}\} = 0$ for $K = K^\star + 1, \ldots, K_{\max}$ in order to conclude the study of $\widehat{K}_n^{\mathrm{G}}$. Minor changes in the previous lines yield the result. □

## APPENDIX C: PROOFS OF EFFICIENCY: UNDERESTIMATION

**C.1. Proof of Theorem 7.** Theorem 7 is first proven under assumption (ii). The modification of the proof under assumption (ii)′ is sketched at the end of this subsection. Let us begin with some useful lemmas.

LEMMA C.1. *Under the assumptions of Theorem 7, the sets $\Lambda_{\alpha,K}$ and $\Gamma_{\alpha,K}$ are measurable and closed in $\mathcal{Q}$ for every $\alpha > 0$ and $K < K^\star$.*

PROOF. The measurability issue is obvious. Set $\alpha > 0$ and $K < K^\star$. We shall actually prove that $\Lambda_{\alpha,K}^c$ is an open set (the same proof applies to $\Gamma_{\alpha,K}$, up to minor changes). To this end, let us point out that the topology on $\mathcal{Q}$ is generated by the collection of open sets

$$\mathcal{O}(f,x,\varepsilon) = \{Q \in \mathcal{Q} : |Qf - x| < \varepsilon\} \qquad (\text{any } f \in \mathcal{L}_\tau, x \in \mathbb{R}, \varepsilon > 0).$$

Choose $Q_0 \in \Lambda_{\alpha,K}^c$, $\alpha', \varepsilon > 0$ such that $\alpha' - 6\varepsilon > \alpha$ and $\sup_{\theta \in \Theta_K} Q_0 \ell_\theta - \sup_{\theta \in \Theta_{K+1}} Q_0 \ell_\theta < -\alpha'$. Let us denote by $T_K$ (resp. $T_{K+1}$) the finite sieve subset of $\Theta_K$ (resp. $\Theta_{K+1}$) for $Q = Q_0$, $\varepsilon$ and $K$ (resp. $K+1$) in assumption (ii). Let us then define the open neighborhood $V$ of $Q_0$ by

$$V = \bigcap_{t \in T_K} \{Q \in \mathcal{Q} : |Q\ell_t - Q_0 \ell_t| < \varepsilon\} \cap \bigcap_{t \in T_{K+1}} \{Q \in \mathcal{Q} : |Q\ell_t - Q_0 \ell_t| < \varepsilon\}.$$

Straightforwardly, whenever $Q \in V$,

$$\sup_{\theta \in \Theta_K} Q\ell_\theta \leq \sup_{\theta \in \Theta_K} Q_0 \ell_\theta + 3\varepsilon$$

and

$$\sup_{\theta \in \Theta_{K+1}} Q_0 \ell_\theta \leq \sup_{\theta \in \Theta_K} Q\ell_\theta + 3\varepsilon,$$

hence, $Q \in \Lambda_{\alpha,K}^c$. So $V$ is an open neighborhood of $Q_0$ included in $\Lambda_{\alpha,K}^c$. This completes the proof of the lemma since $Q_0$ was arbitrarily chosen in $\Lambda_{\alpha,K}^c$. □

LEMMA C.2. *Under the assumptions of Theorem 7, for every $K < K^\star$, $\Pi_K \subset \Lambda_{0,K} \cap \Gamma_{0,K}$ and $P^\star \notin \Lambda_{0,K} \cup \Gamma_{0,K}$.*

PROOF. Let $\tau^*$ be the convex-conjugate of $\tau$, given by $\tau^*(t) = (1 + |t|)\log(1 + |t|) - |t|$ (all $t \in \mathbb{R}$). One can substitute $\tau^*$ for $\tau$ in the definitions (1) of $\mathcal{L}_\tau(P^\star)$ and (3) of $\|\cdot\|_\tau$, yielding a Banach space $(\mathcal{L}_{\tau^*}(P^\star), \|\cdot\|_{\tau^*})$.

30 A. CHAMBAZNow, according to (2.2) in [32], $P^\star|fg| \leq 2\|f\|_\tau \|g\|_{\tau^*}$ [all $f \in \mathcal{L}_\tau(P^\star)$, $g \in \mathcal{L}_{\tau^*}(P^\star)$], so that $\mathcal{L}_{\tau^*}(P^\star)$ can be identified with a subspace of $\mathcal{L}'_\tau(P^\star)$.

Furthermore, it is readily seen that the density of any $P_\theta \in \Pi_{K^\star-1}$ with respect to $P^\star$ belongs to $\mathcal{L}_{\tau^*}(P^\star)$, hence, $\Pi_{K^\star-1} \subset \mathcal{Q}$.

Let us choose $K < K^\star$ and $P_{\theta_0} \in \Pi_K$. Since $p_{\theta_0} \ell_{\theta_0} \in L^1(\mu)$, $\sup_{\theta \in \Theta_{K'}} P_{\theta_0} \ell_\theta = -\inf_{\theta \in \Theta_{K'}} H(P_{\theta_0}|P_\theta) + P_{\theta_0} \ell_{\theta_0} = P_{\theta_0} \ell_{\theta_0}$ when $K' \geq K$. Straightforwardly, $\Pi_K \subset \Lambda_{0,K} \cap \Gamma_{0,K}$.

Besides, $P^\star \in \Lambda_{0,K}$ would yield

$$\sup_{\theta \in \Theta_K} P^\star \ell_\theta - \sup_{\theta \in \Theta_{K+1}} P^\star \ell_\theta = -H(P^\star|\Pi_K) + H(P^\star|\Pi_{K+1}) = 0$$

and $P^\star \in \Gamma_{0,K}$ would yield, in turn,

$$0 \leq \sup_{\theta \in \Theta_K} P^\star \ell_\theta - \sup_{\theta \in \Theta_{K^\star}} P^\star \ell_\theta = -H(P^\star|\Pi_K),$$

where the right-hand side term is negative because $H(P^\star|\cdot)$ achieves its infimum on the compact set $\Pi_K$ and $P^\star \notin \Pi_K$. This completes the proof of the lemma. $\square$

The proof of Theorem 7 follows.

Because $P^\star\{\widehat{K}_n^{\mathrm{L}} < K^\star\} = \sum_{K < K^\star} P^\star\{\widehat{K}_n^{\mathrm{L}} = K\}$, Lemma 1.2.15 of [17] ensures that

$$\limsup_{n \to \infty} n^{-1} \log P^\star\{\widehat{K}_n^{\mathrm{L}} < K^\star\} = \sup_{K < K^\star} \limsup_{n \to \infty} n^{-1} \log P^\star\{\widehat{K}_n^{\mathrm{L}} = K\}.$$

Thus, it suffices to choose $K < K^\star$ and show that

(C.1) $$\limsup_{n \to \infty} n^{-1} \log P^\star\{\widehat{K}_n^{\mathrm{L}} = K\} \leq -I(\Lambda_{0,K}) < 0$$

in order to get (10). Now, for any $\alpha > 0$,

$$\limsup_{n \to \infty} n^{-1} \log P^\star\{\widehat{K}_n^{\mathrm{L}} = K\} \leq \limsup_{n \to \infty} n^{-1} \log P^\star\{\mathbb{P}_n \in \Lambda_{\alpha,K}\} \leq -I(\Lambda_{\alpha,K})$$

by virtue of Theorem 1 and Lemma C.1 $[\mathrm{cl}(\Lambda_{\alpha,K}) = \Lambda_{\alpha,K}]$.

Furthermore, $\{I(\Lambda_{\alpha,K})\}$ nondecreases as $\alpha \downarrow 0$ and it is bounded by $H(\Pi_K|P^\star)$ by virtue of Lemma C.2. Let us denote $L = \lim_{\alpha \downarrow 0} I(\Lambda_{\alpha,K}) \leq I(\Lambda_{0,K}) \leq H(\Pi_K|P^\star)$.

Since $I$ is lower semicontinuous with compact level sets, it achieves its infimum on the closed sets $\Lambda_{\alpha,K}$: let $Q_p \in \Lambda_{1/p,K}$ be such that $I(Q_p) = I(\Lambda_{1/p,K})$ for every $p \geq 1$. For any $q \geq 1$, the set $\mathrm{cl}(\{Q_p : p \geq q\})$ is compact [it is closed in the compact set $\Lambda_{1/q,K} \cap \{Q \in \mathcal{Q} : I(Q) \leq L\}$]. By virtue of the Borel–Lebesgue property, the intersection of the nonincreasing sequence of nonvoid compact sets $\{\mathrm{cl}(\{Q_p : p \geq q\})\}_{q \geq 1}$ is nonvoid too, so $\overline{Q}$ can be chosen in the intersection.

Now, it is readily seen that both $\overline{Q} \in \Lambda_{0,K}$ and $I(\overline{Q}) = I(\Lambda_{0,K}) = L$. Finally, Lemmas 2 and C.2 guarantee that $I(\overline{Q}) > 0$ and yield (C.1), hence, (10).

The proof of (11) for $\widehat{K}_n^{\mathrm{G}}$ is almost identical and is omitted.

*Proof under assumption* (ii)′. Let us assume that $P^\star$-a.s., for all $n \geq 1$, $Z_1, \ldots, Z_n$ are mutually distinct. Then $P^\star$-a.s., $\mathbb{P}_n \in \mathcal{P}'$, where $\mathcal{P}'$ is the subset of $\mathcal{P}$ when adding the condition that $z_1, \ldots, z_p$ must be mutually distinct in the definition of $\mathcal{P}$. Besides, since $I$ is infinite on $\mathcal{P}$, one can substitute $\mathcal{P}'$ for $\mathcal{P}$ in the definition of $\mathcal{Q}$ (see Lemma 4.1.5 in [17]).

The framework introduced for the large deviations principle was intentionally somewhat too simple (for sake of legibility). Under assumption (ii), this is just a matter of convention. When dealing with assumption (ii)′, we must be more careful.

Now, rigorously, $Q \in \mathcal{P}$ is a linear form on $L_\tau(P^\star)$, which has the same definition as $\mathcal{L}_\tau(P^\star)$ except that $P^\star$-almost everywhere equal functions *are not* identified. The topology on $\mathcal{Q}$ is the coarsest one that makes the linear forms $Q \mapsto Qf$ continuous for all $f \in L_\tau(P^\star)$. This change has no effect on $\mathcal{Q} \cap \mathcal{L}'_\tau(P^\star)$. It nevertheless allows to prove that each $Q \in \mathcal{P}'$ is its own open neighborhood in $\mathcal{Q}$.

Indeed, choose $Q_0 = p^{-1} \sum_{i=1}^p \delta_{z_i} \in \mathcal{P}'$. Let $u > 1$ be such that $u/(u-1) < (p+1)/p$ and $V = \bigcap_{i=1}^m \{Q \in \mathcal{Q} : |Q\mathbb{1}\{z_i\} - 1/p| < (up)^{-1}\}$. Then

- $Q_0 \in V$ and $V$ is open.
- If $Q \in V$, then $Q \in \mathcal{P}$ (otherwise, $Q\mathbb{1}\{z_1\} = 0$).
- If $Q = m^{-1} \sum_{i=1}^m \delta_{\zeta_i}$, then $\{\zeta_1, \ldots, \zeta_m\} \supset \{z_1, \ldots, z_p\}$, hence, particularly $m \geq p$.
- Finally, $Q \in V$ yields $|1/m - 1/p| < (up)^{-1}$, which implies, in turn, $m < p+1$, hence, $m = p$ and $Q = Q_0$.

This property allows us to adapt straightforwardly the proof of Lemma C.1 under assumption (ii)′, proving thus the last statement of Theorem 7. □

**C.2. Proof of Theorem 8.** Let us first state some preliminary lemmas.

LEMMA C.3. *Under assumptions* (i) *of Theorem* 7 *and* (iii) *of Theorem* 8, *if* $Q \in \mathcal{Q} \cap \mathcal{L}'_\tau(P^\star)$, *then the function* $\theta \mapsto Q\ell_\theta$ *mapping* $\Theta_{K^\star}$ *to* $\mathbb{R}$ *is continuous over* $\Theta_{K^\star}$.

LEMMA C.4. *Let us choose* $Q \in \mathcal{Q}$ *and* $K \leq K^\star$. *Under assumptions* (i) *of Theorem* 7 *and* (iv) *of Theorem* 8, *the function* $\theta \mapsto Q\ell_\theta$ *mapping* $\Theta_K$ *to* $\mathbb{R}$ *is differentiable on the interior of* $\Theta_K$, *with derivative* $\theta \mapsto Q\dot{\ell}_\theta$.

In order to show Lemma C.3, it is sufficient to prove that $\|\ell_{\theta_p} - \ell_{\theta_0}\|_\tau \to 0$ when $\theta_p \to \theta_0$ in $\Theta_{K^\star}$ (dominated convergence theorem). Lemma C.4 simply relies on the positivity of $Q \in \mathcal{Q}$. Combining both lemmas yields the following:

LEMMA C.5. *Let $Q \in \mathcal{L}'_\tau(P^\star)$ be $P^\star$-singular (i.e., $Q = Q^s$). Then, under the assumptions of Theorem 8, $\theta \mapsto Q\ell_\theta$ is constant over $\Theta_{K^\star}$.*

Consequently, by applying Lemma C.5 to $Q = \overline{Q}$ (see the end of the proof of Theorem 7 in Section C.1), $\overline{Q} \in \Lambda_{0,K}$ yields $\overline{Q}^a \in \Lambda_{0,K} \cap M_1(\mathcal{Z})$. Forwardly,

$$H(\Lambda_{0,K} \cap M_1(\mathcal{Z})|P^\star) \leq H(\overline{Q}^a|P^\star) = I(\overline{Q}^a)$$
$$\leq I(\overline{Q}) = I(\Lambda_{0,K})$$
$$\leq I(\Lambda_{0,K} \cap M_1(\mathcal{Z})) = H(\Lambda_{0,K} \cap M_1(\mathcal{Z})|P^\star),$$

which concludes the proof of Theorem 8 for $\widehat{K}_n^{\mathrm{L}}$. The study of $\widehat{K}_n^{\mathrm{G}}$ goes along the same lines, up to minor changes. □

**C.3. Proof of Theorem 9.** The proof relies heavily on Lemma 5, which is shown at the end of this section. If one resumes the proof of Theorem 8, it is clear that the following proposition straightforwardly yields the result of Theorem 9:

PROPOSITION C.1. *If $\overline{Q} \in \mathcal{L}'_\tau(P^\star) \cap M_1(\mathcal{Z})$ satisfies*

$$H(\overline{Q}|P^\star) = H(\Gamma_{0,K^\star-1} \cap M_1(\mathcal{Z})|P^\star) < \infty,$$

*then, under the assumptions of Theorem 9, $\overline{Q} \in \Pi_{K^\star-1}$ and $H(\overline{Q}|P^\star) = H(\Pi_{K^\star-1}|P^\star)$.*

REMARK C.1. A simple modification of the proof below implies that, under the assumptions of Theorem 9 and for $K = K^\star - 1$,

$$H(\Lambda_{0,K} \cap M_1(\mathcal{Z})|P^\star) = H(\Pi_K|P^\star).$$

The proof cannot be adapted anymore when $K < K^\star - 1$ (the unadaptable argument is pointed out).

PROOF OF PROPOSITION C.1. Let us set $K = K^\star - 1$ and $\overline{Q}$ as described in Proposition C.1. Let us suppose that the assumptions of Theorem 9 are valid. The hard part is to show that $\overline{Q} \in \Pi_K$ since Lemma C.2 guarantees that $\Pi_K \subset \Gamma_{0,K} \cap M_1(\mathcal{Z})$.

Because $\Pi_K$ is compact and $H(\overline{Q}|\cdot)$ is lower semicontinuous, there exists $\overline{P} \in \Pi_K$ (whose density $\mathrm{d}\overline{P}/\mathrm{d}\mu$ is denoted by $\bar{p}$) such that $H(\overline{Q}|\overline{P}) =$ $H(\overline{Q}|\Pi_K)$. According to the definition (14), $\overline{P}$ is the reversed-$H$-projection of $\overline{Q}$ on $\Pi_K$.

We prove hereafter that $\overline{Q} = \overline{P} \in \Pi_K$, which is the expected result.

Let us introduce the subset $\mathcal{C}$ of $M_1(\mathcal{Z}) \cap \mathcal{L}'_\tau(P^\star)$ defined by

$$\mathcal{C} = \{Q : H(Q|P^\star) < \infty\} \cap \left\{Q : Q\log\bar{p} = \sup_{\theta \in \Theta_{K^\star}} Q\ell_\theta\right\}$$

$$\cap \{Q : Q \ll \mu, dQ/d\mu = q, Q\log q < \infty\}.$$

The following properties hold (their simple proofs are omitted):

- $\mathcal{C}$ is convex, $\overline{Q} \in \mathcal{C}$ and $\mathcal{C} \subset \Gamma_{0,K}$.
- $H(\overline{Q}|P^\star) = H(\mathcal{C}|P^\star)$.
- For every $Q \in \mathcal{C}$, $H(Q|\overline{P}) = H(Q|\Pi_K) = H(Q|\Pi_{K^\star})$.

Accordingly,

- $\overline{Q}$ is the $H$-projection of $P^\star$ on $\mathcal{C}$.
- $\overline{P}$ is the reversed-$H$-projection on $\Pi_{K^\star}$ of every $Q \in \mathcal{C}$.

Since the assumptions of Lemma 5 are satisfied, for every $Q \in \mathcal{C}$,

$$H(Q|P^\star) \geq H(Q|\overline{P}) + H(\overline{P}|P^\star)$$

[just choose $P^\star$ in (15)—we point out that this argument is not adaptable when dealing with $\Lambda_{0,K}$ for $K < K^\star - 1$]. Consequently, the characterization of Lemma 4 guarantees that, necessarily, $\overline{P}$ is the $H$-projection of $P^\star$ on $\mathcal{C}$, that is, $\overline{Q} = \overline{P}$, hence, $\overline{Q} \in \Pi_K$. This completes the proof of Proposition C.1. $\square$

PROOF OF LEMMA 5. Since $H(Q|\Pi_{K^\star}) < \infty$, there exists $P_\theta$ such that $H(Q|P_\theta) < \infty$, hence, $Q \ll P_\theta$ and $Q \ll \mu$.

Obviously, if (15) holds, then $H(Q|P) = H(Q|\Pi_{K^\star})$ and $H(Q|P_\theta) = H(Q|P)$ yields $H(P|P_\theta) = 0$, hence, $P_\theta = P$.

Conversely, the exponential nature of the model is needed:

$$p_\theta(z) = h(z)\exp[\theta^T t(z) - \phi(\theta)] \qquad (\text{all } z \in \mathcal{Z}),$$

where $t = (t_1, \ldots, t_{K^\star})$ is a known function on $\mathcal{Z} \subset \mathbb{R}^q$ equipped with Lebesgue measure $\mu$ on Borel sets, $h \in \mathbb{R}^{\mathcal{Z}}$ is measurable, and $\Theta_{K^\star}$ is a convex subset of the convex and open natural parameter space $\Theta = \{\theta \in \mathbb{R}^{K^\star} : \mu(h\exp(\theta^T t)) < \infty\}$; $\phi(\theta) = \log\mu(h\exp(\theta^T t))$ (all $\theta \in \Theta$). Let us emphasize that $\phi$ is convex and differentiable on $\Theta$, with $\dot{\phi}(\theta) = P_\theta t$.

Let $Q$ be chosen as described in the lemma. Let $P$ be its reversed-$H$-projection on $\Pi_{K^\star}$ (with density $dP/d\mu$ denoted by $p$). Inequality (15) is obvious if $H(Q|P_\theta) = \infty$. Consequently, only the parameters $\theta \in \Theta_c = \{\theta \in \Theta : H(Q|P_\theta) < \infty\}$ have to be considered.

Now, it is readily seen that the set $\Theta_c$ is convex. Moreover, it is an open set. Indeed, $Q \log q$ and $Q\ell_\theta$ are finite [because $Q \in \mathcal{L}'_\tau(P^\star)$ and $\ell_\theta \in \mathcal{L}_\tau(P^\star)$], so the decomposition $H(Q|P_\theta) = Q \log q - Q\ell_\theta$ is valid. Assumption (iv) guarantees then that $\Theta_c$ is open.

Besides, because $Q \log q$ and $H(Q|P)$ are finite, (15) is equivalent to $(Q - P) \log p/p_\theta \geq 0$ (any $\theta \in \Theta_c$). Denoting $P$ by $P_{\bar\theta}$ finally implies that (15) is equivalent to

$$(\text{C.2}) \qquad (\bar\theta - \theta)^T (Q - P) t \geq 0 \qquad (\text{all } \theta \in \Theta_c).$$

This concludes that part of the proof.

Now, the decomposition $H(Q|P_\theta) = Q\log q - Q\ell_\theta$ and $H(Q|P) = H(Q|\Pi_{K^\star})$ also imply that

$$(\text{C.3}) \qquad 0 \leq Q \log \frac{p}{p_\theta} < \infty \qquad (\text{all } \theta \in \Theta_c).$$

Let us define $f$ on $\Theta_c$ by $f(\theta) = Q \log p/p_\theta = (\bar\theta - \theta)^T Q t + \phi(\theta) - \phi(\bar\theta)$. Then the convexity of $\phi$ and (C.3) imply that $f$ is a proper convex function on $\Theta_c$. Furthermore, $f$ is differentiable at $\bar\theta$ with gradient $\dot f(\bar\theta) = (P - Q)t$. Since $f$ achieves its minimum at $\bar\theta$ by virtue of (C.3), Theorem 27.4 of [39] applies, hence,

$$(\theta - \bar\theta)^T (-\dot f(\bar\theta)) = (\theta - \bar\theta)^T (Q - P) t \leq 0 \qquad (\text{all } \theta \in \Theta_c).$$

This is exactly (C.2), so the proof is complete. $\square$

## APPENDIX D: PROOFS OF EFFICIENCY: OVERESTIMATION

Let us denote by $\Delta_K = D(K+1) - D(K) > 0$ and $K = K^\star + 1$,

$$(\text{D.1}) \quad P^\star\{\widehat K_n^{\text{L}} > K^\star\} \leq P^\star\left\{\sup_{\theta \in \Theta_{K^\star+1}} \mathbb{P}_n \ell_\theta - \sup_{\theta \in \Theta_{K^\star}} \mathbb{P}_n \ell_\theta \geq n^{-1} v_n \Delta_{K^\star}\right\}$$

$$(\text{D.2}) \qquad \leq P^\star\left\{\sup_{g \in \mathcal{G}_K^{\text{a}}} |(\mathbb{P}_n - P^\star)g| \geq n^{-1} v_n \Delta_{K^\star}\right\},$$

by virtue of (A.1). Also, the peeling device inequality (A.2) of the same proposition implies that expression given by (D.1) can be bounded by

$$(\text{D.3}) \qquad P^\star\left\{\left(\sup_{g \in \mathcal{G}_K^{\text{b}}} |(\mathbb{P}_n - P^\star)g|\right)^2 \geq n^{-1} v_n \Delta_{K^\star}\right\}.$$

In the rest of this paper, we shall focus on (16) in Theorem 10 [on the basis of the overestimation probability upper bound (D.2)]. The proof of (18) in Theorem 11 [on the basis of the overestimation probability upper bound (D.3)] is similar and is omitted.

Let us define $\Lambda^\infty = \{b \in \ell^\infty(\mathcal{G}_K^a) : \|b\|_{\mathcal{G}_K^a} \geq \Delta_{K^\star}\}$. It is closed for the uniform topology on $\ell^\infty(\mathcal{G}_K^a)$. Since the assumptions of Theorem 2 are satisfied,

$$\limsup_{n\to\infty} nv_n^{-2} \log P^\star\{\widehat{K}_n^L > K^\star\} \leq \limsup_{n\to\infty} nv_n^{-2} \log P^\star\{(nv_n^{-1})(\mathbb{P}_n - P^\star)^\infty \in \Lambda^\infty\}$$
$$\leq -\inf\{J(b) : b \in \Lambda^\infty\}.$$

Let us prove that the right-hand side term above is negative.

Suppose indeed, on the contrary, that the infimum is zero: this implies $0 \in \Lambda^\infty$, which is obviously not true. If the infimum were zero, then there would exist a sequence $\{b_p\}$ of elements of $\ell^\infty(\mathcal{G}_K^a)$ such that $b_p \in \Lambda^\infty$ and $J(b_p) \leq 1/p$. Consequently, there would exist a sequence $\{Q_p\}$ of elements of $M(\mathcal{Z})$ such that, for every $p \geq 1$, $Q_p \ll P^\star$ (with derivative $dQ_p/dP^\star$ denoted by $q_p$) and both $P^\star q_p^2/2 \leq J(b_p) + 1/p \leq 2/p$ and $Q_p^\infty = b_p$. Thus, for any $g \in \mathcal{G}_K^a$,

$$(b_p g)^2 = (P^\star q_p g)^2 \leq (P^\star q_p^2)(P^\star g^2) \leq (4/p)\left(\sup_{g \in \mathcal{G}_K^a} P^\star g^2\right)$$

by virtue of the Cauchy–Schwarz inequality. Now, $\mathcal{G}_K^a$ is $P^\star$-Donsker, hence, it is totally bounded in $L^2(P^\star)$, and the above display implies that $\|b_p\|_{\mathcal{G}_K^a} = o(1)$. Consequently, $0 \in \Lambda^\infty$ as a limit of a sequence of elements of the closed set $\Lambda^\infty$.

This completes the proof of (16) of Theorem 10.

The proof of (17) in Theorem 10 [which parallels the proof of (19) in Theorem 11] is very similar. Once again, the union bound and Lemma 1.2.15 of [17] imply that

$$\limsup_{n\to\infty} nv_n^{-2} \log P^\star\{\widehat{K}_n^G > K^\star\}$$
$$= \sup_K \limsup_{n\to\infty} nv_n^{-2} \log P^\star\{\widehat{K}_n^G = K\}$$
$$\leq \sup_K \limsup_{n\to\infty} nv_n^{-2} \log P^\star\left\{\sup_{\theta \in \Theta_K} \mathbb{P}_n \ell_\theta - \sup_{\theta \in \Theta_{K^\star}} \mathbb{P}_n \ell_\theta \geq n^{-1} v_n \Delta_K\right\}$$

($\sup_K$ stands for $\sup_{K^\star < K \leq K_{\max}}$). This bound is handled as the bound (D.1) above, hence, the final result. $\square$

## APPENDIX E: PROOFS FOR THE BENCHMARK EXAMPLES

**E.1. Proof of Lemma 6.** Lemma E.1 allows us to focus on the restrictions of $\ell_\theta$'s to a well-chosen compact set of $\mathcal{Z}$.

LEMMA E.1. *Let $\psi \in \mathbb{R}_+^{\mathbb{R}_+}$ be an increasing nonnegative function such that $\psi(x)/x \to \infty$ as $x \to \infty$. Let us assume that $(u - l)$ is continuous on* $\mathcal{Z} \subset \mathbb{R}^q$, $(u-l)(z) \to \infty$ as $|z| \to \infty$, and that $\psi(u-l) \in \mathcal{L}_\tau(P^\star)$. Then $(u-l) \in \mathcal{M}_\tau(P^\star)$ and, for all $\varepsilon > 0$ and $Q \in \mathcal{Q}$, there exists a compact subset $C$ of $\mathcal{Z}$ such that $Q(u-l)\mathbb{1}\{C^c\} < \varepsilon$.

PROOF. It is easily verified that $(u-l) \in \mathcal{M}_\tau(P^\star)$. Besides, the set $\{z \in \mathcal{Z} : (u-l)(z) \leq M\}$ is compact for any $M > 0$ and $Q(u-l)\mathbb{1}\{(u-l) > M\} \leq \frac{M}{\psi(M)} Q\psi(u-l) \to 0$ as $M \to \infty$. $\square$

Of course, the assumptions of Lemma E.1 are satisfied in the LM example [here, $\psi(x) = x^{1+c}$ (any $x \geq 0$)]. Thus, let us set $K \geq 1$, $Q \in \mathcal{Q}$ and $\varepsilon > 0$. There exists a compact set $C$ of $\mathcal{Z}$ such that, for every $\theta, t \in \Theta_K$,

$$(\text{E.1}) \quad |Q(\ell_\theta - \ell_t)| \leq Q|\ell_\theta - \ell_t|\mathbb{1}\{C\} + Q(u-l)\mathbb{1}\{C^c\} \leq Q|\ell_\theta - \ell_t|\mathbb{1}\{C\} + \varepsilon.$$

Now, Ascoli's theorem ensures that $\{\ell_\theta \mathbb{1}\{C\} : \theta \in \Theta_K\}$ is precompact in the set of the continuous functions on $C$ equipped with the uniform norm. Consequently, there exists a finite subset $T$ of $\Theta_K$ such that, for every $\theta \in \Theta_K$, there exists $t \in T$ such that $\sup_{z \in C} |\ell_\theta(z) - \ell_t(z)| \leq \varepsilon$. Straightforwardly, for any $\theta \in \Theta_K$, there exists $t \in T$ such that the left-hand side term of (E.1) is bounded by $2\varepsilon$. This completes the proof. $\square$

**E.2. Proof of Lemma 7.** Let us suppose, on the contrary, that

$$(\text{E.2}) \quad H(P^\star|\Pi_K) \leq H(P^\star|\Pi_{K+1}),$$

that is, that equality holds. Lower semicontinuity of $H(P^\star|\cdot)$ and compactness of $\Pi_K$ ensure the existence of $P_0 = P_{\theta_0} \in \Pi_K$ such that $H(P^\star|P_0) = H(P^\star|\Pi_K)$. Let us denote $f_0(x) = f_{\theta_0}(x) = \sum_{k=1}^K m_k \mathbb{1}\{x \in \tau_k\}$ (all $x \in \mathcal{X}$). Now, equality (20) and $\|f^\star - f_0\|^2 = \sum_{k=1}^K P(f^\star - m_k)^2 \mathbb{1}\{\tau_k\}$ imply that $m_k = Pf^\star \mathbb{1}\{\tau_k\}/P(\tau_k)$ for $k = 1, \ldots, K$. Let us prove that $f^\star = f_0$, hence, $P^\star \in \Pi_K$.

Indeed, (E.2) ensures that, for any $1 \leq k_0 \leq K$, for any subset $S$ of $\tau_{k_0}$ with positive $P$-measure,

$$Pf^{\star 2} - \sum_{k=1}^K \frac{(Pf^\star \mathbb{1}\{\tau_k\})^2}{P(\tau_k)} \leq Pf^{\star 2} - \sum_{1 \leq k \neq k_0 \leq K} \frac{(Pf^\star \mathbb{1}\{\tau_k\})^2}{P(\tau_k)}$$
$$- \left( \frac{(Pf^\star \mathbb{1}\{S\})^2}{P(S)} + \frac{(Pf^\star \mathbb{1}\{\tau_{k_0} \setminus S\})^2}{P(\tau_{k_0} \setminus S)} \right)$$

or, equivalently,

$$\frac{(Pf^\star \mathbb{1}\{S\})^2}{P(S)} + \frac{(Pf^\star \mathbb{1}\{\tau_{k_0} \setminus S\})^2}{P(\tau_{k_0} \setminus S)} \leq \frac{(Pf^\star \mathbb{1}\{\tau_{k_0}\})^2}{P(\tau_{k_0})}.$$

Thus, first expansion of the right-hand side term and then factorization yield

$$\frac{P(\tau_{k_0} \setminus S)}{P(S)}(Pf^\star 1\{S\})^2 + \frac{P(S)}{P(\tau_{k_0} \setminus S)}(Pf^\star 1\{\tau_{k_0} \setminus S\})^2$$
(E.3)
$$\leq 2(Pf^\star 1\{S\})(Pf^\star 1\{\tau_{k_0} \setminus S\}).$$

Now, the basic inequality $2ab \leq (au)^2 + (bu^{-1})^2$ (all $a,b \in \mathbb{R}$ and positive $u$) together with (E.3) ensure [take $u^2 = P(\tau_{k_0} \setminus S)/P(S)$] that equality holds in (E.3). Consequently, for any subset $S$ of $\tau_{k_0}$ with positive $P$-measure,

$$\frac{Pf^\star 1\{S\}}{P(S)} = \frac{Pf^\star 1\{\tau_{k_0} \setminus S\}}{P(\tau_{k_0} \setminus S)} = \frac{Pf^\star 1\{\tau_{k_0}\} - Pf^\star 1\{S\}}{P(\tau_{k_0 \setminus S})},$$

hence, for any subset $S$ of $\tau_{k_0}$,

$$Pf^\star 1\{S\} = \frac{P(S)}{P(\tau_{k_0})} Pf^\star 1\{\tau_{k_0}\}.$$

The choice $S = S_+ = \{x \in \tau_{k_0} : f^\star(x) > Pf^\star 1\{\tau_{k_0}\}/P(\tau_{k_0})\}$ yields $P(S_+) = 0$. The choice $S = S_- = \{x \in \tau_{k_0} : f^\star(x) < Pf^\star 1\{\tau_{k_0}\}/P(\tau_{k_0})\}$ yields, in turn, $P(S_-) = 0$, hence, finally $P(S_0) = P(\tau_{k_0})$, where $S_0 = \{x \in \tau_{k_0} : f^\star(x) = Pf^\star 1\{\tau_{k_0}\}/P(\tau_{k_0})\}$ (i.e., $f^\star$ $P$-a.s. constant on $\tau_{k_0}$). This concludes the proof because $k_0$ is arbitrary. □

**E.3. Proof of Lemma 8.** Let us set $K \geq 1$, $Q \in \mathcal{Q} \cap \mathcal{L}'_\tau(P^\star)$ (with decomposition $Q = Q^a + Q^s$ according to Lemma 1) and $\varepsilon > 0$. Because $Q^a \ll P^\star$, there exists $\delta > 0$ such that, for any measurable $F$, $P^\star(F) \leq \delta$ yields $Q^a(F) \leq \varepsilon$.

Now, it was emphasized in Section 5.2 that $(u-l)^{1+c} \in \mathcal{L}_\tau(P^\star)$ for some $c > 0$, hence, Lemma E.1 applies with $\psi(x) = x^{1+c}$ (all $x \geq 0$). So, there exists a compact set $C$ of $\mathcal{Z}$ such that, for every $\theta, t \in \Theta_K$,

$$|Q(\ell_\theta - \ell_t)| \leq Q|\ell_\theta - \ell_t|1\{C\} + Q(u-l)1\{C^c\}$$
$$\leq Q|\ell_\theta - \ell_t|1\{C\} + \varepsilon$$
(E.4)
$$= Q^a|\ell_\theta - \ell_t|1\{C\} + \varepsilon$$
$$\leq MQ^a|f_\theta - f_t| + \varepsilon,$$

where the equality holds because $(\ell_\theta - \ell_t)1\{C\}$ is bounded and $M$ is a constant which depends only on $l, u$ (via $C$) and $\mathcal{M}$.

Furthermore, the Borel–Lebesgue property of compact sets guarantees that there exists a finite subset $T$ of $\Theta_K$ such that the union over $t \in T$ of the balls of center $t$ and radius $\delta$ covers $\Theta_K$. Let us set $t \in T$ [$t = (\tau^0, m^0)$] and $\theta \in \Theta_K$ [$\theta = (\tau^1, m^1)$] with $d(t, \theta) \leq \delta$. It can be assumed without loss of generality that $P(\tau_j^0 \Delta \tau_j^1) \leq \delta$ and $|m_j^0 - m_j^1| \leq \delta$ for all $j = 1, \ldots, K$. Consequently, with notation, $M' = \sup\{|m| : m \in \mathcal{M}\}$, for any $x \in \mathcal{X}$,

$$|f_\theta - f_t|(x) \leq \sum_{j=1}^{K} |m_j^0 - m_j^1| \mathbb{1}\{x \in \tau_j^0 \cap \tau_j^1\} + M'(K-1) \sum_{j=1}^{K} \mathbb{1}\{x \in \tau_j^0 \Delta \tau_j^1\}$$

$$\leq K\delta + M'(K-1) \sum_{j=1}^{K} \mathbb{1}\{x \in \tau_j^0 \Delta \tau_j^1\},$$

hence,

$$Q^a |f_\theta - f_t| \leq K\delta + M'(K-1) \sum_{j=1}^{K} Q^a(\tau_j^0 \Delta \tau_j^1 \times \mathbb{R}).$$

Besides, $P^\star(\tau_j^0 \Delta \tau_j^1 \times \mathbb{R}) = P(\tau_j^0 \Delta \tau_j^1) \leq \delta$ finally yields $Q^a(\tau_j^0 \Delta \tau_j^1 \times \mathbb{R}) \leq \varepsilon$. By invoking (E.4), $|Q(\ell_\theta - \ell_t)| \leq M''\varepsilon$, for a constant $M''$ depending only on $K$, $l, u$ (via $C$) and $\mathcal{M}$. This completes the proof. $\square$

**Acknowledgments.** This work was done while I was affiliated with the University Paris-Sud and France Télécom Recherche & Développement. I wish to express my gratitude to my Ph.D. advisor Elisabeth Gassiat. I would also like to thank Stéphane Boucheron, Raphaël Cerf, Christian Léonard, Pascal Massart and Jamal Najim for helpful discussions. I am especially grateful to one of the referees for his suggestions and careful reading.

MAP5 CNRS UMR 8145
UNIVERSITÉ RENÉ DESCARTES
45 RUE DES SAINTS-PÈRES
75270 PARIS CEDEX 06
FRANCE
E-MAIL: chambaz@univ-paris5.fr